\documentclass[12pt]{article}

\usepackage{fullpage}
\usepackage{verbatim}
\usepackage{calrsfs}
\usepackage{amsthm}
\usepackage{amsbsy}
\usepackage{amssymb}
\usepackage{amsfonts}
\usepackage[dvips]{lscape}
\usepackage{hyperref} 
\usepackage{amscd}    

\begin{document}
\title{Big Line Bundles Over Arithmetic Varieties}
\author{Xinyi Yuan}
\maketitle

\theoremstyle{plain}
\newtheorem{thm}{Theorem}[section]
\newtheorem*{conj}{Conjecture}
\newtheorem{cor}[thm]{Corollary}
\newtheorem*{corr}{Corollary}
\newtheorem{lem}[thm]{Lemma}
\newtheorem{pro}[thm]{Proposition}
\newtheorem{definition}[thm]{Definition}
\newtheorem*{thmm}{Theorem}

\theoremstyle{remark} \newtheorem*{remark}{Remark}

\newcommand{\lb}{\mathcal{L}}             
\newcommand{\mb}{\mathcal{M}}
\newcommand{\nb}{\mathcal{N}}
\newcommand{\fb}{\mathcal{F}}
\newcommand{\eb}{\mathcal{E}}
\newcommand{\tb}{\mathcal{T}}
\newcommand{\ob}{\mathcal{O}}
\newcommand{\ib}{\mathcal{I}}
\newcommand{\jb}{\mathcal{J}}

\newcommand{\xb}{\mathcal{X}}             

\newcommand{\lbb}{\overline{\mathcal{L}}}             
\newcommand{\mbb}{\overline{\mathcal{M}}}
\newcommand{\nbb}{\overline{\mathcal{N}}}
\newcommand{\ebb}{\overline{\mathcal{E}}}
\newcommand{\tbb}{\overline{\mathcal{T}}}
\newcommand{\obb}{\overline{\mathcal{O}}}
\newcommand{\lbt}{\widetilde{\mathcal{L}}}
\newcommand{\mbt}{\widetilde{\mathcal{M}}}
\newcommand{\chern}{\hat{c}_1}           
\newcommand{\height}{h_{\lbb}}           
\newcommand{\cheight}{\hat{h}_{\lb}}     

\newcommand{\gs}[1]{\Gamma(X, #1)}        
\newcommand{\rank}{\mathrm{rank}}         
\newcommand{\vol}[1]{\mathrm{vol}(#1)}    
\newcommand{\lnorm}[1]{\|#1\|_{L^2}}      
\newcommand{\supnorm}[1]{\|#1\|_{\mathrm{sup}}}    
\newcommand{\chil}{\chi_{_{L^2}}}
\newcommand{\chisup}{\chi_{\sup}}
\newcommand{\divv}{\mathrm{div}}         
\newcommand{\Spec}{\mathrm{Spec}}         

\newcommand{\xcv}{X_{\CC_v}^{\mathrm{an}}}   
\newcommand{\xkv}{X_{K_v}^{\mathrm{an}}}     

\newcommand{\RR}{\mathbb{R}}
\newcommand{\PP}{\mathbb{P}}      
\newcommand{\QQ}{\mathbb{Q}}
\newcommand{\CC}{\mathbb{C}}
\newcommand{\ZZ}{\mathbb{Z}}      

\newcommand{\pn}{\mathbb{P}^n}

\newcommand{\gm}{\mathbb{G}_{\mathrm{m}}}
\newcommand{\gmn}{\mathbb{G}_{\mathrm{m}}^n}

\tableofcontents

\section{Introduction}
In this paper, we prove a bigness theorem (Theorem \ref{main}) in the setting of Arakelov theory as an
arithmetic analogue of a classical theorem of Siu. It is also an extension of the arithmetic Hilbert-Samuel
formula implied by the arithmetic Riemann-Roch theorem of Gillet-Soul\'e \cite{GS2} and an estimate on analytic torsions of Bismut-Vasserot \cite{BV}. Our treatment of arithmetic bigness is based on the theory of arithmetic ampleness by Zhang \cite{Zh1}.

This bigness result has a lot of consequences in the
equidistribution theory initiated by Szpiro-Ullmo-Zhang \cite{SUZ}. We
will generalize to algebraic dynamics the archimedean
equidistribution by \cite{SUZ}, the non-archimedean equidistribution by
Chambert-Loir \cite{Ch2}, and the equidistribution of small
subvarieties by Baker-Ih \cite{BI} and Autissier \cite{Au2}.

The equidistribution theorem in \cite{SUZ} was proved by a variational principle (cf. \cite{Ch3}), where the key is to use the arithmetic Hilbert-Samuel formula to produce small sections. The formula works under the assumption that the curvature of the line bundle giving the polarization is strictly positive, since any small perturbation of the line bundle still have positive curvature. Such an assumption is also necessary in \cite{Ch2} if the variety has dimension greater than one, while the complete result in the case of curves is obtained there by using a result of Autissier \cite{Au1} which we will recall later.

However, in algebraic dynamics (e.g. multiplicative groups), the curvature is usually only semipositive and even a small perturbation may result in a somewhere negative curvature. Then the arithmetic Hilbert-Samuel is invalid in this case. Our bigness theorem solves this problem, since it works for negative curvatures.

Our proof of the bigness theorem follows a strategy similar to the one used to prove the arithmetic Hilbert-Samuel formula by Abbes and Bouche \cite{AB}. The crucial analytic part is the estimate of the distortion function of $N\lbb-j\mbb$ in Proposition \ref{3.3}. It is implied by its ample case (Theorem \ref{3.2}) proved by Bouche \cite{Bo} and Tian \cite{Ti}.

\subsection{Equidistribution over Algebraic Dynamics}

\subsubsection*{Projective Spaces}
Let $K$ be a number field, and $\overline K$ be the algebraic closure of $K$. Fix an embedding $\overline K \rightarrow \CC$. Let $\pn$ be the projective space over $K$, and $\phi:\pn\rightarrow \pn$ be an endomorphism with coordinate $\phi=(f_0, f_1, \cdots, f_n)$, where $f_0, f_1, \cdots, f_n$ are homogeneous polynomials of degree $q>1$ without non-trivial common zeros.

For any algebraic point $x=(z_0, z_1, \cdots, z_n) \in \PP^n(\overline K)$, the \textit{naive height} of $x$ is
$$h_{\rm naive}(x):=\frac{1}{[L:K]}\sum_v \log \max\{|z_0|_v, |z_1|_v,\cdots, |z_m|_v \},$$
where $L$ is a finite extension of $K$ containing all the coordinates $z_0, z_1, \cdots, z_m$, and the summation is over all normalized valuations $|\cdot|_v$ of $L$.

The \textit{canonical height} with respect to $\phi$ is defined by Tate's limit
$$h_{\phi}(x)=\lim_{N\rightarrow \infty} \frac{1}{q^N}h_{\rm naive}(\phi^N(x)).$$
One can show that the limit always exists.

The canonical height has the following nice property: $h_{\phi}(x)\geq 0$ and $h_{\phi}(x)=0$ if and only if $x$ if \textit{preperiodic}. Here we say a point is preperiodic if its orbit $\{x, \phi(x), \phi^2(x), \cdots\}$ is finite.

To state our equidistribution theorem, we make some simple definitions related to sequences of algebraic points of
$\pn(\overline K)$.
\begin{enumerate}
\item A sequence $\{x_m\}_{m\geq 1}$ of algebraic points is \textit{small} if $h_{\phi}(x_m)\rightarrow 0$ as       $m\rightarrow \infty.$
\item A sequence $\{x_m\}_{m\geq 1}$ of algebraic points is \textit{generic} if no infinite subsequence of               $\{x_m\}$ is contained in a proper closed subvariety of $\pn$.
\item Let $\{x_m\}_{m\geq 1}$ be a sequence of algebraic points and $d\mu$ a \textit{probability measure} over the complex manifold $\pn(\CC)$, i.e., a measure of total volume one. We say that the Galois orbits of $\{x_m\}$ are \textit{equidistributed} with respect to $d\mu$ if the probability measure $\displaystyle\mu_{x_m}:=\frac{1}{\#O(x_m)}\sum_{x\in O(x_m)}\delta_x$ converges weakly to $d\mu$ over $\pn(\CC)$, where $O(x_m)$ is the orbit of $x_m$ under the Galois group $\mathrm{Gal}(\overline K/K)$, and $\delta_x$ is the Dirac measure at $x\in \pn(\CC).$
\end{enumerate}

We can also define \textit{the canonical probability measure} $d\mu_{\phi}$ over $\pn(\CC)$ by Tate's limit. It is a probability measure that satisfies $\phi^*d\mu_{\phi}=q^{\dim(X)}d\mu_{\phi}$ and $\phi_*d\mu_{\phi}=d\mu_{\phi}$, which determine $d\mu_{\phi}$ uniquely.

The following theorem is a special case of Theorem \ref{5.7} in this paper.
\begin{thmm}
Suppose $\{x_m\}_{m\geq 1}$ is an infinite sequence of algebraic points in $X$ which is generic and small. Then the Galois orbits of $\{x_m\}$ are equidistributed with respect to the canonical probability measure $d\mu_{\phi}$ over $\pn(\CC)$.
\end{thmm}

\subsubsection*{Generalities}
We actually prove the equidistribution for any algebraic dynamical systems in Theorem \ref{5.7}. For a complete introduction of algebraic dynamics and related equidistribution we refer to \cite{Zh5}.

Let $K$ be a number field. An \textit{algebraic dynamical system} over $K$ is a projective variety $X$ over $K$ endowed with an endomorphism $\phi:X\rightarrow X$ which satisfies a polarization condition making it like the polynomial map over $\pn$. By Tate's limit, we have the same notion of \textit{canonical height} and \textit{canonical probability measure}. And thus we have the dynamical equidistribution over $X$. See Section \ref{dynamics} for more details.

Now we are going to consider three special cases:
\begin{enumerate}
\item Abelian varieties. When $X$ is an abelian variety and $\phi=[2]$ is multiplication by 2, we get a dynamical system. A point is preperiodic if and only if it is torsion. The canonical height is exactly the Neron-Tate height, and the canonical probability measure is exactly the probability Haar measure over the complex torus $X(\CC)$. Our equidistribution in this case is exactly the one in \cite{SUZ}, which was crucial in the proof of the Bogomolov conjecture by Ullmo \cite{Ul} and Zhang \cite{Zh3}. See \cite{Zh4} for an abstract of this subject.

\item Multiplicative groups. When $X=\gmn$ and $\phi=[2]$, we get a dynamical system over multiplicative groups. To compactify it, embed $\gmn$ in $\pn$ by the natural way and extend $\phi$ to a dynamics over $\pn$. Actually  $\phi:(z_0, z_1, \cdots, z_n)\mapsto(z_0^2, z_1^2, \cdots, z_n^2)$. The curvature is semipositive here, which can't be handled by \cite{SUZ}. That is why the proof of the Bogomolov conjecture by Zhang \cite{Zh1} and the proof of equidistribution by Bilu \cite{Bi} were independent of each other and could not follow the idea of Ullmo and Zhang. However, our new result puts this case in the framework of Ullmo and Zhang.

\item Almost split semi-abelian varieties. In \cite{Ch1}, Chambert-Loir
proved equidistribution and Bogomolov conjecture over almost split
semi-abelian varieties. The equidistribution was proved by
choosing certain nice perturbation of the line bundle which
preserves the semipositivity of the curvature. As in the
multiplicative case, it can be handled by our uniform treatment.
\end{enumerate}

\subsection{A Generic Equidistribution Theorem}
The above equidistribution theorem is implied by the following generic equidistribution theorem in Arakelov geometry.
The basic references for Arakelov geometry are \cite{Ar}, \cite{Fa}, \cite{GS1} and \cite{Zh1}.

Let $K$ be a number field, $X$ be a projective variety of dimension $n-1$ over $K$, and $\lb$ be a line bundle over $X$.
Fix an embedding $\overline K\rightarrow \CC_v$ for each place $v$, where $\CC_v$ is the completion of the algebraic closure of $K_v$.

We use the language of \textit{adelic metrized line bundles} by Zhang \cite{Zh1, Zh2}. Recall that an \textit{adelic metric} over $\lb$ is a $\CC_v$-norm $\|\cdot\|_v$ over the fibre $\lb_{\CC_v}(x)$ of each algebraic point $x\in X(\overline K)$  for each place $v$ of $K$ satisfying certain continuity and coherence conditions.

All the metrics we consider are induced by models or uniform limits of metrics induced by models. Suppose $(\xb,\lbt)$ is an $O_K$-model of $(X, \lb^e)$, i.e., $\xb$ is an integral scheme projective and flat over $O_K$ and $\lbt$ is a Hermitian line bundle over $\xb$ such that the generic fibre of $(\xb,\lbt)$ gives $(X, \lb^e)$. For any non-archimedean place $v$, a point $x\in X(\overline K)$  extends to $\tilde x: \mathrm{Spec}(O_{\CC_v})\rightarrow \xb_{O_{\CC_v}}$. Then
$\displaystyle (\tilde x^* \lbt_{O_{\CC_v}})^{\frac{1}{e}}$
gives a lattice in $\lb_{\CC_v}(x)$, which induces a $\CC_v$-norm and thus an adelic metric. Such a metric is called an \textit{algebraic metric}. It is called \textit{semipositive} if $\lbt$ has semipositive curvatures at all archimedean places and non-negative degree on any complete vertical curve of $\xb$.

An adelic metric over $\lb$ is \textit{semipositive} if it is the uniform limit of some sequence of semipositive algebraic metrics over $\lb$.

\theoremstyle{plain}\newtheorem*{thm 5.1}{Theorem \ref{5.1}}
\begin{thm 5.1}[Equidistribution of Small Points]
Suppose $X$ is a projective variety of dimension $n-1$ over a number field $K$, and $\lbb$ is a metrized line bundle over $X$ such that $\lb$ is ample and the metric is semipositive. Let $\{x_m\}$ be an infinite sequence of
algebraic points in $X(\overline K)$ which is generic and small. Then for any place $v$ of $K$, the Galois orbits of the sequence $\{x_m\}$ are equidistributed in the analytic space $\xcv$ with respect to the probability measure $d\mu_v=c_1(\lbb)_v^{n-1}/\deg_{\lb}(X)$.
\end{thm 5.1}

We explain several terms in the theorem:
\begin{enumerate}
\item The definitions of a generic sequence and equidistribution
are the same as before. 
\item Using the semipositive line bundle $\lbb$, one can define the height $\height(Y)$ of any closed
subvariety $Y$ of $X$. Namely, 
$$\height(Y)=\frac{\chern(\lbb)^{\dim Y+1}|_{\overline Y}}{(\dim Y+1)\deg_{\lb}(\overline Y)},$$
where $\overline Y$ is the closure of $Y$ in the scheme $X$. A sequence $\{x_m\}$ of algebraic points in $X(\overline K)$ is called \textit{small} if $\height(x_m)\rightarrow \height(X)$. 
\item For archimedean $v$,
the space $\xcv$ is the corresponding complex analytic space, and
the measure $c_1(\lbb)_v^{n-1}$ is essentially the volume form induced by
the hermitian metric of $\lbb$ at $v$. See \cite[Proposition 3.1.5]{Zh5} for example. In this case, our theorem generalizes \cite[Theorem 3.1]{SUZ}. 

\item For non-archimedean $v$, the theorem
generalizes the recent work of Chambert-Loir \cite{Ch2}. Here $\xcv$ is
the Berkovich space (cf. \cite{Be}) associated to the variety
$X_{\CC_v}$.  Chambert-Loir constructs the $v$-\textit{adic
canonical measure} $c_1(\lbb)_v^{n-1}$ and generalizes the
equidistribution of \cite{SUZ} to the $v$-adic case.  We follow
Chambert-Loir's notion of canonical measures. 
\item Another ingredient in our non-archimedean treatment is a theorem of Gubler
\cite{Gu} that any continuous real-valued function over $\xcv$ can be
approximated by \textit{model functions} induced by certain formal
models. In our case, all model functions are induced by arithmetic
varieties, which puts the problem in the framework of Arakelov
theory. Finally, we obtain a proof analogous to the archimedean
case, in which model functions play the role of smooth functions.
\end{enumerate}

\begin{remark}
The results in \cite{SUZ} and \cite{Ch2} assume the strict positivity of the metric at the place where equidistribution is considered in the case that $\dim X>1$.
\end{remark}

\subsection{Arithmetic Bigness}
Our bigness theorem is the key to deal with negative curvatures. Here we state it and explain how it works.

\subsubsection*{Siu's Theorem}
Let $X$ be a projective variety of dimension $n$ defined over a field, and $\lb$ be a line bundle over $X$. If $\lb$ is ample, then when $N$ is large enough, the Hilbert function $h^0(\lb^{\otimes N})=\dim \Gamma(X, \lb^{\otimes N})$ is a polynomial in $N$ of degree $n$. Notice that ampleness is stable under pull-back via finite morphisms, but not via birational morphisms.

Another useful notion for line bundles is bigness, which is stable under pull-back via dominant generically finite morphisms. The line bundle $\lb$ is \textit{big} if and only if there exists a constant $c>0$ such that
$h^0(\lb^{\otimes N}) > cN^n$ for all $N$ large enough. See \cite{La} for more details of bigness.

Denote by $c_1(\lb_1)\cdots c_1(\lb_n)$ the intersection number of the line bundles $\lb_1, \cdots, \lb_n$ over $X$. The following is a basic theorem of Siu \cite{Si}. See also \cite[Theorem 2.2.15]{La}.

\begin{thmm}[Siu]
Let $\lb$, $\mb$ and $\eb$ be three line bundles over a projective variety $X$ of dimension $n$. Assume that $\lb$ and $\mb$ are ample. Then
$$h^0(\eb+N(\lb-\mb))
\geq \frac{c_1(\lb)^n-n\cdot c_1(\lb)^{n-1} c_1(\mb)}{n!}N^n +O(N^{n-1}).$$
In particular, $\lb-\mb$ is big if $c_1(\lb)^n> n \cdot c_1(\lb)^{n-1}c_1(\mb)$.
\end{thmm}

Here we write tensor product of line bundles additively, like the case of divisors. For example, $\eb+N(\lb-\mb)$ means $\eb\otimes(\lb\otimes \mb^{\otimes (-1)})^{\otimes N}$.

\subsubsection*{Arithmetic Bigness}
One arithmetic analogue of the classical $h^0$ is $\chi_{\sup}$ (cf. \cite{GS2}).
See also Section \ref{notation} for an explanation. Our direct analogue of Siu's
theorem gives a nice expansion of $\chi_{\sup}$. Its accuracy
allows it to play the role of the arithmetic Hilbert-Samuel
formula in equidistribution.

Let $X$ be an arithmetic variety of dimension $n$, and let $\lbb$ be a hermitian line bundle over $X$. We say that $\lbb$ is \textit{strongly big} if there exists a constant $c>0$ such that $\chi_{\sup}(\lbb^{\otimes N})> cN^n$ for all $N$ large enough.

Note that there is a nice arithmetic theory of \textit{ample line bundles} by Zhang \cite{Zh1}. Namely, a hermitian line bundle $\lbb$ is \textit{ample} if the following three conditions are satisfied:

(a) $\lb_\QQ$ is ample in the classical sense;

(b) $\lbb$ is \textit{relatively semipositive}: the curvature of $\lbb$ is semipositive and $\deg(\lb|_C) \geq 0$ for any closed curve $C$ on any special fibre of $X$ over $\mathrm{Spec}(\mathbb{Z})$;

(c) $\lbb$ is \textit{horizontally positive}: the intersection number $\chern(\lbb|_Y)^{\dim Y}>0$ for any horizontal irreducible closed subvariety $Y$.

Now we have the following main theorem which has the same appearance as Siu's theorem:
\theoremstyle{plain}\newtheorem*{thm 1.2}{Theorem \ref{main}}
\begin{thm 1.2}[Main Theorem]
Let $\lbb$, $\mbb$ and $\ebb$ be three hermitian line bundles over an arithmetic variety $X$ of dimension $n$. Assume that $\lbb$ and $\mbb$ are ample. Then
$$\chi_{\sup}(\ebb+N(\lbb-\mbb))
\geq \frac{\chern(\lbb)^n-n\cdot\chern(\lbb)^{n-1}\chern(\mbb)}{n!}N^n +o(N^n).$$
In particular, $\lbb-\mbb$ is strongly big if $\chern(\lbb)^n > n\cdot\chern(\lbb)^{n-1}\chern(\mbb)$.
\end{thm 1.2}

We can compare this theorem with the arithmetic Hilbert-Samuel formula. Actually the former is like a bigness version of the latter.

\begin{thmm}[Arithmetic Hilbert-Samuel]
Let $\lbb$ and $\ebb$ be two line bundles over an arithmetic
variety $X$ of dimension $n$. If $\lbb$ is relatively semipositive
and $\lb_{\QQ}$ is ample, then
\begin{equation}
\chi_{\sup}(\ebb+N\lbb)=\frac{\chern(\lbb)^n}{n!}N^n+o(N^n) \ , \ N\rightarrow \infty.
\end{equation}
\end{thmm}

The Hilbert-Samuel formula was originally proved in Gillet-Soul\'e \cite{GS2} by combining an estimate of Bismut-Vasserot \cite{BV}. The above one is an extension by Zhang \cite{Zh1}. The original one was also proved by Abbes-Bouche \cite{AB} using a more straight-forward method. We will extend the method in this paper to prove Theorem \ref{main}.

Now let us see how the bigness theorem works in proving the equidistribution. The variational principle in \cite{SUZ} is to consider the bundle $\lbb(\epsilon f)=(\lb, e^{-\epsilon f}\|\cdot\|_{\lbb})$, the same line bundle $\lb$ with metric multiplied by $e^{-\epsilon f}$ at $v$. Here $f$ is any smooth function over the analytic space $\xcv$, and $\epsilon>0$ is a small number.

The strategy is to write $\ob(f)=\mbb_1-\mbb_2$ for ample hermitian line bundles $\mbb_1$ and $\mbb_2$, where $\ob(f)$ is the trivial line bundle with metric $\|1\|=e^{-f}$. Then $\lbb(\epsilon f)=(\lbb+\epsilon\mbb_1)-\epsilon\mbb_2$ is a difference of two ample hermitian line bundles and we can apply Theorem \ref{main} to this difference. Note that $\epsilon \mbb_2$ is small, and the leading term given by the theorem actually approximates $\chern(\lbb(\epsilon f))^n$ up to an error $O(\epsilon^2)$.

\subsection{Structure of this Paper}
The structure of this paper is as follows. In Section 2 we state and prove the main theorem (Theorem \ref{main}), and explore several basic properties of arithmetic bigness. 

Sections 2.3-2.5 give a proof of the main theorem. The outline of the proof is clear in Section 2.5. Some preliminary results on arithmetic volumes (resp. analytic estimate) are proved in Section 2.3 (resp. Section 2.4). We also reduce the problem to certain good case in Section \ref{reduction}.

Section 3 gives a detailed treatment of the equidistribution theory, which reveals the importance of Theorem \ref{main}. The readers that are more interested in algebraic dynamics may assume Theorem \ref{main} and jump directly to Section 3. 

\subsection*{Acknowledgments}
I am very grateful to my advisor Shou-wu Zhang who introduced this
subject to me. I am also indebted to him for lots of helpful
conversations and encouragement during the preparation of this
paper. I would like to thank Brian Conrad for clarification of
many concepts in rigid analytic geometry and Aise Johan de Jong 
for his help with algebraic geometry. Finally, I would like to
thank Zuoliang Hou, Xander Faber, Ming-lun Hsieh and Qi Li for
their useful discussions.

During the revision of this paper, Moriwaki \cite{Mo3} proved the
continuity of the volume functions of hermitian line bundles using
the technique of the proof of Theorem \ref{main} and Theorem
\ref{1.3} in this paper.

\section{Arithmetic Bigness}

\subsection{Notations and Conventions}\label{notation}

By an \textit{arithmetic variety} $X$ of dimension $n$, we mean an integral scheme $X$, projective and flat over $\mathrm{Spec}(\mathbb{Z})$ of absolute dimension $n$. We say that $X$ is \textit{generically smooth} if the generic fibre $X_{\mathbb{Q}}$ is smooth. In this paper we don't assume $X$ to be generically smooth, and use generic resolution to relate the general case to the generically smooth case by Hironaka's theorem. See \cite{Zh1} for more on generic resolutions. In Proposition \ref{2.2}, one will see that resolution of singularities preserves bigness of line bundles very well, so we actually don't need to worry about singularities on the generic fibre.

A \textit{metrized line bundle} $\lbb=(\lb, |\cdot|)$ over $X$ is an invertible sheaf $\mathcal{L}$ over $X$ together with a hermitian metric $|\cdot|$ on each fibre of $\mathcal{L}_{\mathbb{C}}$ over $X_\CC$. We say this metric is \textit{smooth} if the pull-back metric over $f^*\lb$ under any analytic map $f: \{z\in \CC^{n-1}:|z|<1 \} \rightarrow X_\CC$ is smooth in the usual sense. We call $\lbb$ a \textit{hermitian line bundle} if its metric is smooth and invariant under complex conjugation. For a hermitian line bundle $\lbb$, we say the metric or the curvature of $\lbb$ is \textit{semipositive} if the curvature of $f^*\lbb$ with the pull-back metric under any analytic map $f: \{z\in \CC^{n-1}:|z|<1 \} \rightarrow X_\CC$ is semipositive definite.

Let $\lbb_1, \lbb_2, \cdots, \lbb_n$ be $n$ hermitian line bundles over $X$. Choose any generic resolution $\pi:\widetilde{X} \rightarrow X$. Then the \textit{intersection number} $\hat{c}_1({\lbb_1})\hat{c}_1({\lbb_2})\cdots \hat{c}_1({\lbb_n})$ is defined to be $\hat{c}_1(\pi^*{\lbb_1})\hat{c}_1(\pi^*{\lbb_2})\cdots \hat{c}_1(\pi^*{\lbb_n})$, where the latter is the usual arithmetic intersection number defined in \cite{GS1}. This definition is independent of the choice of the generic resolution (cf. \cite{Zh1}).

For any section $s\in \Gamma(X, \lb)_\RR =\Gamma(X, \lb)\otimes_{\ZZ}\RR \subset \Gamma(X_\CC, \lb_\CC),$ one has the supremum norm $\supnorm{s}=\sup_{z\in X_{\mathbb{C}}}|s(z)|$. Define a basic invariant
$$h^0(\lbb)=\log \#\left\{s \in \Gamma(X, \lb): \supnorm{s}<1\right\}.$$
\noindent Picking any Haar measure on $\Gamma(X, \lb)_\RR$, define the \textit{arithmetic volume}
$$\chi_{\sup}(\lbb)=\log\frac{\vol{B_{\sup}}}{\vol{\gs{\lb}_{\RR}/\gs{\lb}}},$$
\noindent where $B_{\sup}=\{ s\in\gs{\lb}_{\RR}: \supnorm{s}<1  \}$ is the corresponding unit ball. It is easy to see that this definition is independent of the choice of the Haar measure.

Zhang studied arithmetic ampleness in \cite{Zh1}. Recall that a hermitian line bundle $\lbb$ is \textit{ample} if the following three conditions are satisfied:

(a) $\lb_\QQ$ is ample;

(b) $\lbb$ is \textit{relatively semipositive}: the curvature of $\lbb$ is semipositive and $\deg(\lb|_C) \geq 0$ for any closed curve $C$ on any special fibre of $X$ over $\mathrm{Spec}(\mathbb{Z})$;

(c) $\lbb$ is \textit{horizontally positive}: the intersection number $\chern(\lbb|_Y)^{\dim Y}>0$ for any horizontal irreducible closed subvariety $Y$.

Note that the second condition in (b) means $\lb$ is nef over any
special fibre in the classical sense. By Kleiman's theorem, it is
equivalent to $c_1(\lb|_Y)^{\dim Y}\geq 0$ for any vertical
irreducible closed subvariety $Y$. See \cite[Theorem 1.4.9]{La}.

The arithmetic Hilbert-Samuel formula is true for ample line bundles, and thus we can produce a lot of small sections. Zhang proved an arithmetic Nakai-Moishezon theorem on this aspect. The following result is a combination of Theorem 4.2 and Theorem 3.5 in Zhang \cite{Zh1}. 

\begin{thmm}[Zhang]
Let $\lbb$ be an ample hermitian line bundle on an arithmetic variety $X$. Assume that there exists an embedding
\ $i: X_\CC \rightarrow Y$ to a projective mainfold $Y$, and an ample line bundle $\lb'$ on $Y$ with $i^*\lb'=\lb_\CC$, such that the metric of $\lb_\CC$ can be extended to a hermitian metric on $\lb'$ with semipositive curvature.

Then for any hermitian line bundle $\ebb$ over $X$, the $\ZZ$-module $\gs{\eb+N\lb}$ has a basis consisting of \textit{strictly effective sections} for $N$ large enough. 
\end{thmm}

In particular, the assumption is automatic is $X$ is generically smooth. See also Zhang \cite[Corollary 4.8]{Zh1}.

Here an \textit{effective section} is a nonzero section with
supremum norm less than or equal to 1. We call a line bundle
\textit{effective} if it admits an effective section. If the
supremum norm of the section is less than 1, the section and the
line bundle are said to be \textit{strictly effective}.

In the end, we state a fact telling that conditions (a) and (b)
are not far from ampleness. More precisely, if $\lbb=(\lb,
\|\cdot\|)$ is such that $\lb_\QQ$ is ample and $\lbb$ is
relatively semipositive, then the hermitian line bundle
$\lbb(c)=(\lb, \|\cdot\|_c= \|\cdot\|e^{-c})$ is ample for $c$
large enough. In fact, since $\lb_\QQ$ is ample, we can assume
there exist sections $s_1, \cdots, s_r \in \Gamma(X, \lb)$ which
are base-point free over the generic fibre. Fix a $c$ such that
$s_1, \cdots, s_r$ are strictly effective in $\lbb(c)$. Now we
claim that $\lbb(c)$ is ample. We need to show that
$\chern(\lbb(c)|_Y)^{\dim Y}>0$ for any horizontal irreducible
closed subvariety $Y$. Assume $X$ is normal by normalization. We
can find an $s_j$ such that $\mathrm{div}(s_j)$ does not contain
$Y$, and thus
\begin{eqnarray*}
\chern(\lbb(c)|_Y)^{\dim Y}
&=&\chern(\lbb(c)|_{\mathrm{div}(s_j)|_Y})^{\dim
Y-1}-\int_{Y_{\CC}}\log \|s_j\|_c \ c_1(\lbb)^{\dim Y-1}\\
&>&\chern(\lbb(c)|_{\mathrm{div}(s_j)|_Y})^{\dim Y-1}.
\end{eqnarray*}
Now the proof can be finished by induction on $\dim Y$. This fact
is used in Lemma \ref{5.3} when we apply Theorem \ref{main}.

\subsection{Big Line Bundles}
Now we define two notions of arithmetic bigness, which are weaker than ampleness but allow more flexibility.

\begin{definition}\label{1.1}
Let $X$ be an arithmetic variety of dimension $n$, and let $\lbb$ be a hermitian line bundle over $X$.
We say that $\lbb$ is big if there exist a positive integer $N_0$ and a positive number $c$ such that
$$ h^0(N\lbb) > cN^n$$ for any integer $N>N_0$.
We say that $\lbb$ is strongly big if there exist a positive integer $N_0$ and a positive number $c$ such that
$$  \chi_{\sup}(N\lbb)> cN^n$$
for any integer $N>N_0$.
\end{definition}

\begin{remark}
\begin{enumerate}

\item Moriwaki \cite{Mo2} defines that $\lbb$ is big if $\lb_{\QQ}$ is
big in the classical sense and some positive power of $\lbb$ is
strictly effective. It turns out that his definition is equivalent
to ours. See Corollary \ref{1.4} below.

\item Minkowski's theorem gives $h^0(N\lbb)\geq\chi_{\sup}(N\lbb)
+O(N^{n-1}),$ and thus ``strongly big" implies ``big". Its
converse is not true in general. An example will be showed at the
end of this section. \item Either notion of bigness is invariant
under dominant generically finite morphisms; i.e., the pull-back
bundle of a big (resp. strongly big) line bundle via a dominant
generically finite morphism is still big (resp. strongly big).
\item In the two-dimensional case, Autissier \cite[Proposition
3.3.3]{Au1} proved a strong result for general line bundles. Namely,
$\displaystyle\chi_{\sup}(N\lbb)\geq \frac{\chern(\lbb)^2}{2}N^2
+o(N^2)$ for any hermitian line bundle $\lbb$ over an arithmetic
surface such that $\deg(\lb_\QQ)>0.$ It tells us that $\lbb$ is
strongly big if and only if $\deg(\lb_\QQ)>0$ and
$\chern(\lbb)^2>0$ by Corollary \ref{1.4} below.
\end{enumerate}
\end{remark}

The main theorem in this paper is the following:

\begin{thm}\label{main}
Let $\lbb$, $\mbb$ and $\ebb$ be three hermitian line bundles over an arithmetic variety $X$ of dimension $n$. Assume that $\lbb$ and $\mbb$ are ample. Then
$$\chi_{\sup}(\ebb+N(\lbb-\mbb))
\geq \frac{\chern(\lbb)^n-n\cdot\chern(\lbb)^{n-1}\chern(\mbb)}{n!}N^n +o(N^n).$$
In particular, $\lbb-\mbb$ is strongly big if $\chern(\lbb)^n > n\cdot\chern(\lbb)^{n-1}\chern(\mbb)$.
\end{thm}

\begin{remark}
The inequality above has an error term. We explain its meaning here. Suppose $F, G, H$ are real-valued functions defined on the positive integers, and $H$ is positive-valued. Then the equality $\displaystyle F(N) \geq G(N)+o(H(N))$ means that there exists a function $R$ such that 
$\displaystyle F(N) \geq G(N)+R(N) $ and $\displaystyle R(N)=o(H(N)), \ N\rightarrow \infty $. We have similar understanding if we replace $\geq$ by $\leq$, or replace $o(H(N))$ by $O(H(N))$. 
\end{remark}

\

The theorem will be proved in Section \ref{proof}. But now we will state two properties of bigness. In the classical case, one has: big=ample+effective. More precisely, a line bundle is big if and only if it has a positive tensor power isomorphic to the tensor product of an ample line bundle and an effective line bundle. For the details see \cite{La}. In the arithmetic case, we have a similar result.

\begin{thm}\label{1.3}
A hermitian line bundle $\lbb$ is big if and only if $N \lbb=\mbb+\tbb$ for some positive integer $N$, some ample hermitian line bundle $\mbb$ and some effective hermitian line bundle $\tbb$.
\end{thm}

We will prove this theorem in Section \ref{proof} after proving the main theorem. The proof is similar to some part of the proof of our main theorem. The key is to use the Riemann-Roch theorem in \cite{GS3} to relate $h^0$ to $\chi_{\sup}$. The following corollary gives more descriptions of big line bundles. And it also says that arithmetic bigness implies classical bigness over the generic fibre.

\begin{cor}\label{1.4}
Let $\lbb$ be a hermitian line bundle over an arithmetic variety. The following are equivalent:
\begin{enumerate}
\item[(1)] $\lbb$ is big.
\item[(2)] $N \lbb=\mbb+\tbb$ for some positive integer $N$, some ample hermitian line bundle $\mbb$ and some effective hermitian line bundle $\tbb$.
\item[(3)] For any line bundle $\ebb$, the line bundle $N\lbb+\ebb$ is effective when $N$ is large enough.
\item[(4)] $\lb_{\QQ}$ is big over $X_{\mathbb{Q}}$ in the classical sense and $N\lbb$ is strictly effective for some positive integer $N$.
\end{enumerate}
\end{cor}

\begin{proof}
$(1)\Longleftrightarrow(2)$. It is Theorem \ref{1.3}.

$(3)\Longrightarrow(2)$. It is trivial by setting $\ebb=-\mbb$, where $\mbb$ is any ample hermitian line bundle.

$(2)\Longrightarrow(3)$. Suppose $N \lbb=\mbb+\tbb$ as in (2). Then $rN\lbb+\ebb=(\ebb+r\mbb)+r\tbb.$ Because $\mbb$ is ample, $\ebb+r\mbb$ is effective for $r$ large enough, and thus $rN\lbb+\ebb$ is effective for $r$ large enough. Replacing $\ebb$ by $\ebb+k\lbb$ for $k=0, 1, \cdots, N-1$, we see that $N'\lbb+\ebb$ is effective when $N'$ is large enough.

Property (4) is Moriwaki's definition of big line bundles, and $(3)\Leftrightarrow(4)$ is Proposition 2.2 in \cite{Mo2}. For convenience of readers, we still include it here.

$(2)\Longrightarrow(4)$. Assume $N \lbb=\mbb+\tbb$ as in (2). It is easy to see that $rN\lbb$ is strictly effective for some integer $r>0$. By $N \lb_{\QQ}=\mb_{\QQ}+\tb_{\QQ}$, we see $\lb_{\QQ}$ is big by classical theory.

$(4)\Longrightarrow(2)$. Assume that there exists a section $s\in \gs{N\lb}$ with $\supnorm{s}<1$. Since $\lb_{\QQ}$ is big, the line bundle $-\mb_{\QQ}+N'\lb_{\QQ}$ is effective for some integer $N'>0.$ It follows that $-\mb+N'\lb$ has a regular section $t$. Now $\supnorm{s^rt}\leq \supnorm{s}^r\supnorm{t}<1$ for $r$ large enough. That means $-\mbb+(N'+rN)\lbb$ is effective.
\end{proof}

\begin{remark}
Arithmetically big line bundles share many properties with the
classical big line bundles. An important one is the continuity of
the volume function
$$\widehat{\mathrm{vol}}(\lbb)=\limsup_{N\rightarrow\infty}\frac{ h^0(N \lbb)}{N^n/n!},$$
which is proved by Moriwaki in the recent work \cite{Mo3}. His proof
follows the same strategy as our proof of Theorem \ref{main}
and Theorem \ref{1.3} here.
\end{remark}

To end this section, we give an example that a line bundle is big
but not strongly big. Suppose $X=\PP^1_\ZZ=\mathrm{Proj} \
\ZZ[x_0, x_1]$ and $\tb=O(1)$. Pick a constant $0<c<e^{-1}$, and
define a metric over $\tb$ by
$$\|s(x_0, x_1)\|=\frac{|s(x_0, x_1)|}{\sqrt{|x_0|^2+c|x_1|^2}},$$
where $s(x_0, x_1)$ is considered as a homogeneous linear polynomial in $x_0$ and $x_1$. It is easy to see that the metric is well defined. And the section $s_0(x_0, x_1)=x_0$ is effective.

Let $\displaystyle z=\frac{x_1}{x_0}$ be the usual affine coordinate on $X-V(x_0)$. Direct computation shows that the curvature form
$\displaystyle c_1(\tbb)=\frac{ic}{2\pi} \frac{\mathrm{d}z\wedge \mathrm{d}\bar{z}}{(1+c|z|^2)^2}$ is positive and $\displaystyle \chern(\tbb)^2=\frac{1}{2}(1+\log c)<0$.

Let $\mbb$ be any ample hermitian line bundle over $X$. For $m>0$,
the line bundle $\lbb=\mbb+m\tbb$ is big (ample+effective) and satisfies the arithmetic Hilbert-Samuel formula. But when $m$ is large enough, the leading coefficient $\displaystyle\frac{1}{2}\chern(\lbb)^2=\frac{1}{2}(\chern(\mbb)+m\chern(\tbb))^2$ in the arithmetic Hilbert-Samuel formula is negative. We conclude that $\lbb$ is not strongly big.

\subsection{Arithmetic Volumes}\label{volume}

In this section, we consider general normed modules and list their basic properties in Proposition \ref{2.1} which will be an important tool to read volume information from exact sequences. An important example in this class is the supremum norm and the $L^2$-norm for sections of a hermitian line bundle.

As the first application, we show that strong bigness over an arithmetic variety is implied by strong bigness over its generic resolution in Proposition \ref{2.2}. By this, we reduce the problem to generically smooth arithmetic varieties in Section \ref{reduction}.

\subsubsection{Normed Modules}
By a \textit{normed $\ZZ$-module} $M$ we mean a finitely generated $\ZZ$-module $M$ together with an $\RR$-norm $\|\cdot\|$ on $M_{\RR}=M\otimes{\RR}$. For such an $M$, define
$$h^0(M, \|\cdot\|)=\log \# \{m\in M: \|m\| < 1 \}.$$
\noindent Denote by $M_\mathrm{tor}$ the torsion part of $M$, and by $M_\mathrm{free}$ the free part of $M$. Identify $M_\mathrm{free}$ with the image of $M$ in $M_{\RR}$. Then $M_\mathrm{free}$ is naturally a full lattice in $M_{\RR}$.
Define
$$\chi(M, \|\cdot\|)=\log\frac{\vol{B(M)}}{\vol{M_{\RR}/M_\mathrm{free}}}+\log\#M_\mathrm{tor},$$
where $B(M)=\{m\in M_{\RR}: \|m\| < 1 \}$ is the unit ball for the norm. Define $\chi(M)$ to be $\log \#M_\mathrm{tor}$ if $M$ is torsion. Note that $\chi(M)$ does not depend on the Haar measure chosen over $M_{\RR}$. 
Sometimes, we omit the dependence on the metric and simply write $h^0(M)$ and $\chi(M)$ if no confusion occurs.

The norm associated to $M$ is \textit{quadratic} if it is an inner product on $M_{\RR}$. In this case we call $M$ a \textit{quadratically normed $\ZZ$-module}. If $m_1, m_2, \cdots, m_r$ is a $\ZZ$-basis of $M_{\mathrm{free}}$, then
$$\displaystyle\chi(M)=\log\frac{V(r)}{\sqrt{\det(\left\langle m_j, m_k\right\rangle)_{1\leq j,k\leq r}}}+\log \#M_\mathrm{tor},$$
where 
$$\displaystyle V(r):=\pi^{\frac{r}{2}}/\Gamma(\frac{r}{2}+1)$$
is the volume of the unit ball in the Euclidean space $\RR^r$. Stirling's formula implies that 
$$\log V(r)=-r\log r+O(r), \quad r\rightarrow \infty.$$
We will use this result to control error terms coming from the results below.

Let $M$ be any normed $\ZZ$-module. Then $\chi(M)$ and $h^0(M)$ are related by a Riemann-Roch theorem of Gillet-Soul\'e \cite{GS3}. Fix a $\ZZ$-basis $m_1, \cdots, m_r$ of $M_\mathrm{free}$. This basis identifies $M_\RR$ with $\RR^r$. Then $B(M)=\{m\in M_{\RR}: \|m\| < 1 \}$ is a convex symmetric body in $\RR^r.$ Define
$$h^1(M)=h^1(M,\|\cdot\|)=\log \#\left\{(a_1,\cdots, a_r)\in \ZZ^r: 
\left|\sum_{i=1}^r a_ib_i\right|<1, \forall\ m=\sum_{i=1}^r b_im_i\in B(M)\right\}.$$
One can check that it is independent of the choice of the $\ZZ$-basis $m_1, \cdots, m_r$.
The result of Gillet-Soul\'e is as follows:
\begin{thm} \cite[Theorem 1]{GS3}\label{lattice rr}
When the normed module $M$ varies, 
$$h^0(M)-h^1(M)=\chi(M)+O(r\log r), \quad r\rightarrow \infty.$$
Here $r$ is the rank of $M$, and the error term $O(r\log r)$ depends only on $r$.
\end{thm}

Now we consider some basic properties of $\chi$ and $h^0$, and some easy consequences of the above theorem.
\begin{pro}\label{2.1}
In the following, denote $r=\rank (M),\ r'=\rank (M'),\ r''=\rank (M'')$.
\begin{itemize}
\item[(1)]
For any normed module $M$, we have $h^0(M) \geq  \chi(M)-r\log 2$.

\item[(2)]
If a finitely generated $\ZZ$-module $M$ has two norms $\|\cdot\|_1$ and $\|\cdot\|_2$ such that $\|\cdot\|_1 \leq \|\cdot\|_2$, then 
$$\chi(M, \|\cdot\|_1)\geq \chi(M, \|\cdot\|_2),  \quad \quad h^0(M, \|\cdot\|_1)\geq h^0(M, \|\cdot\|_2),$$
and
$$0\leq h^0(M, \|\cdot\|_1) -h^0(M, \|\cdot\|_2) \leq  \chi(M, \|\cdot\|_1)-\chi(M, \|\cdot\|_2)+O(r\log r).$$

\item[(3)]
If a finitely generated $\ZZ$-module $M$ has two norms $\|\cdot\|_1$ and $\|\cdot\|_2$ such that 
$\|\cdot\|_1 = \alpha\|\cdot\|_2$ for some $\alpha >0$, then
\begin{itemize}
\item[(a)] $\displaystyle \ \chi(M, \|\cdot\|_1) = \chi(M, \|\cdot\|_2) -r\log \alpha $,
\item[(b)] $\displaystyle h^0(M, \|\cdot\|_1)  = h^0(M, \|\cdot\|_2) -r\log \alpha +O(r\log r)$.
\end{itemize}

\item[(4)]
Let $0 \rightarrow M' \rightarrow M \rightarrow M'' \rightarrow 0 $ be an exact sequence of normed modules, i.e., the sequence is exact as $\ZZ$-modules, and the norms on $M'_\RR$ and on $M''_\RR$ are respectively the subspace norm and quotient norm induced from $M_\RR$. 
\begin{itemize}
\item[(a)]
If furthermore the norms are quadratic, then
$$ \chi(M)-\chi(M')-\chi(M'')=\log V(r)-\log V(r')-\log V(r'').$$
In particular, one has $\chi(M) \leq \chi(M')+\chi(M'')$.
\item[(b)]
For general norms, we have
$$h^0(M)\leq h^0(M')+h^0(M'')+O(r'\log r'),$$
where the error term $O(r'\log r')$ depends only on $r'$.
\end{itemize}

\item[(5)] If a quadratically normed module $M$ can be generated by elements with norms not greater than a positive constant $c$, then $\chi(M)\geq \log V(r)-r\log c$.

\item[(6)] Let $f: M'\rightarrow M$ be an injection of quadratically normed $\ZZ$-modules that is norm-contractive, i.e., $\|f(m')\| \leq \|m'\|$ for all $m'\in M'_\RR$. If $M$ can be generated by elements with norms not greater than a positive constant $c$, then 
\begin{itemize}
\item[(a)] $\displaystyle \ \chi(M')\leq \chi(M)- \log V(r)+\log V(r')+(r-r')\log c$,
\item[(b)] $\displaystyle  h^0(M')\leq h^0(M)$.
\end{itemize}

\end{itemize}
\end{pro}

\begin{proof}
\begin{itemize}
\item[(1)]
It is just Minkowski's theorem.

\item[(2)]
We only need to show
$$h^0(M, \|\cdot\|_1) -h^0(M, \|\cdot\|_2) \leq  \chi(M, \|\cdot\|_1)-\chi(M, \|\cdot\|_2)+O(r\log r).$$
It is implied by Theorem \ref{lattice rr} and the fact that $h^1(M,\|\cdot\|_1)\leq h^1(M,\|\cdot\|_2)$.

\item[(3)]
Equality in (a) follows from its definition and (b) is implied by (a) using the last result of (2). See also \cite[Proposition 4]{GS3}.

\item[(4)]
\begin{itemize}
\item[(a)]
$M''_\RR$ is isomorphic to the orthogonal complement of $M'_\RR$ in $M_\RR$ with the induced subspace norm. The result follows from the fact that $\log(\vol{M_{\RR}/M_\mathrm{free}}/\#M_\mathrm{tor})$ is additive if the volume elements are induced by the norms.

It remains to check $\log V(r)\leq \log V(r')+\log V(r'')$. Instead of using the formula for $V(r)$, we propose a geometric way. Let 
$$B(r)=\{(x_1,x_2,\cdots, x_r)\in \RR^r: x_1^2+\cdots+x_r^2 \leq 1\}$$
be the unit ball in the Euclidean space $\RR^r$ for any $r$. Then the map 
$$(x_1,x_2,\cdots, x_r)\mapsto ((x_1,x_2,\cdots, x_{r'}),(x_{r'+1},x_{r'+2},\cdots, x_r) )$$
gives an embedding $B(r) \hookrightarrow B(r')\times B(r'')$. The map keeps the volume elements, so
$\vol{B(r)} \leq \vol{B(r')}\vol{B(r'')}$. It gives the inequality.

\item[(b)]
Denote $L(M)=\{m\in M: \|m\| < 1 \}$, and similarly for $L(M')$ and $L(M'')$. Consider the induced map $p: L(M)\rightarrow L(M'')$. For any $y\in L(M'')$ and $x\in p^{-1}(y)$, the set $p^{-1}(y)-x=\{z-x: z\in p^{-1}(y)\}$ is contained in $L_2(M'):=\{m\in M': \|m\| < 2 \}$. This gives $\#p^{-1}(y)\leq \#L_2(M')$, and thus $\#L(M)\leq (\#L_2(M'))\cdot(\#L(M''))$. Take logarithm and use (3) (b).
\end{itemize}

\item[(5)] 
We can assume that $c=1$ and $M$ is torsion free. By the condition we can find $r$ elements $m_1, m_2, \cdots, m_r \in M$ with $\|m_j\|\leq 1$ which form a $\ZZ$-basis of a submodule $M'$ of finite index. 
Since $$\chi(M)\geq\chi(M')=\log\frac{V(r)}{\sqrt{\det(\left\langle m_j, m_k\right\rangle)_{1\leq j,k\leq r}}},$$
it suffices to show $\displaystyle \det(\left\langle m_j, m_k\right\rangle) \leq 1$. The matrix  
$\displaystyle A=(\left\langle m_j, m_k\right\rangle)_{1\leq j,k\leq r}$ is symmetric and positive definite, so it has $r$ positive eigenvalues $x_1, x_2, \cdots, x_r$. We have 
$$\displaystyle \sum_{j=1}^r x_j=tr(A)=\sum_{j=1}^r \left\langle m_j, m_j\right\rangle \leq r.$$
By the arithmetic-geometric mean inequality, $\displaystyle\det(A)=\prod_{j=1}^r x_j \leq 1$.

\item[(6)] 
The inequality (b) is trivial by definition. We put it with (a) here because we need to compare them later.
As for (a), it suffices to show the case that $M'$ is endowed with the induced subspace norm. Let $M''=M/M'$ be endowed with the quotient norm. Apply (4)(a), and apply (5) for $M''$.
 
\end{itemize}
\end{proof}

\begin{remark}
As the referee points out, the inequality $\det(A) \leq 1$ is some special case of Hadamard's inequality.
\end{remark}

\

\subsubsection{Arithmetic Volumes}

Let $X$ be an arithmetic variety of dimension $n$ and $\lbb$ be a line bundle over $X$. Then the supremum norm $\supnorm{\cdot}$ makes $\Gamma(X, \mathcal{L})_{\RR}$ a normed module. Apparently, it is not quadratic. However, one can define an \textit{$L^2$-norm} which is quadratic and closely related to the supremum norm. 
Fix a measure $d\mu$ on $X(\CC)$, which is assumed to be the push-forward measure of a pointwise positive measure on some resolution of singularity of $X(\CC)$. One defines the $L^2$-norm by 
$$\displaystyle\lnorm{s}=\left(\int |s(z)|^2d\mu\right)^{1/2}, \quad s\in\Gamma(X, \mathcal{L})_\CC.$$
Then this $L^2$-norm makes $\Gamma(X, \mathcal{L})$ a quadratically normed module. 

These two norms induces four invariants $h^0_{\sup}(\lbb)$, $h^0_{L^2}(\lbb)$, $\chi_{\sup}(\lbb)$, $\chi_{L^2}(\lbb)$ for $\lbb$. The following theorem says that these two norms are equivalent in some sense.

\begin{cor}\label{result of gromov}
\begin{itemize}
\item[(1)]
If $\lbb$ is ample, then all $h^0_{\sup}(N\lbb)$, $h^0_{L^2}(N\lbb)$, $\chi_{\sup}(N\lbb)$, $\chi_{L^2}(N\lbb)$ have  the same expansion
$$\displaystyle\frac{\chern(\lbb)^n}{n!}N^n+o(N^n).$$
The same is true for $\ebb+N\lbb$ for any hermitian line bundle $\ebb$.

\item[(2)]
Let $\ebb$, $\lbb$, $\mbb$ be three hermitian line bundles. Then
\begin{eqnarray*}
\chi_{\sup}(\ebb+N\lbb-j\mbb)&=& \chi_{L^2}(\ebb+N\lbb-j\mbb)+O(N^{n-1}\log (N+j)),\\
h^0_{\sup}(\ebb+N\lbb-j\mbb) &=& h^0_{L^2}(\ebb+N\lbb-j\mbb)+O(N^{n-1}\log (N+j)).
\end{eqnarray*}
\end{itemize}
\end{cor}
\begin{proof}
\begin{itemize}
\item[(1)]
By the Nakai-Moishezon type theorem of Zhang \cite[Corollary 4.8]{Zh1}, the
$\ZZ$-module $\gs{N\lb}$ has a basis consisting of effective sections for $N$ large enough. It follows that $h^1(\Gamma(X, N\lb))=0$ under both the supremum norm and the $L^2$-norm. Now the result follows from the arithmetic Hilbert-Samuel formula and Theorem \ref{lattice rr} above.

\item[(2)]
The Gromov inequality in Proposition \ref{3.4} asserts that $c(N+j)^{-n}\supnorm{\cdot} \leq \lnorm{\cdot} \leq \supnorm{\cdot}$. Apply Proposition \ref{2.1} (2), (3).
\end{itemize}
\end{proof}

In the end, we prove a result that enables us to replace $X$ by its generic resolution in next subsection.
\begin{pro}\label{2.2}
Let $\lbb$ and $\ebb$ be two hermitian line bundles over an arithmetic variety $X$ of dimension $n$. Let $\pi:\widetilde{X} \rightarrow X$ be any birational morphism from another arithmetic variety $\widetilde X$ to $X$. Then
$$\chisup(\ebb\otimes\lbb^{\otimes N}) \geq \chisup(\pi^*\ebb\otimes\pi^*\lbb^{\otimes N})+ o(N^n).$$
\end{pro}

\begin{proof}
For simplicity, we assume that $\ebb$ is trivial. The general case is proved in the same way with minor work.

Firstly, we can reduce to the case that $\pi:\widetilde{X} \rightarrow X$ is finite. In fact, consider the Stein factorization $\widetilde{X}\stackrel{p}\rightarrow X'\stackrel{\pi'}\rightarrow X$, where $X'=\mathrm{Spec}(\pi_*\ob_{\widetilde{X}})$ is finite over $X$. One has $p_*\ob_{\widetilde{X}}=\ob_{X'}$. For any hermitian line bundle $\tbb$ over $X'$,
$$\Gamma(\widetilde{X}, p^*\tb)= \Gamma(X', p_*(p^*\tb))=\Gamma(X', \tb \otimes p_*\ob_{\widetilde{X}})=\Gamma(X', \tb)$$
by projection formula. The isomorphism $\Gamma(\widetilde{X}, p^*\tb)= \Gamma(X', \tb)$ is actually an isometry under the supremum norms. Therefore, $\chisup(\widetilde{X}, \pi^*\lbb^{\otimes N})= \chisup(X', \pi'^*\lbb^{\otimes N})$. So it suffices to show the same result for the morphism $\pi':X' \rightarrow X$, which is finite.

Secondly, it suffices to prove $\chil(\lbb^{\otimes N}) \geq \chil(\pi^*\lbb^{\otimes N})+ o(N^n)$ by choosing nice measures on $X$ and $\widetilde{X}$. Suppose $\widetilde{X}' \rightarrow \widetilde{X}$ is a generic resolution of $\widetilde{X}$. Fix a pointwise positive measure over $\widetilde{X}'$, which induces push-forward measures over $\widetilde{X}$ and $X$. These measures define $L^2$-norms for line bundles over them. By the above corollary, a bound on $\chil$ is equivalent to the same bound on $\chisup$.

Now assume $\pi:\widetilde{X} \rightarrow X$ is finite, and $\widetilde{X}, X$ are endowed with measures as above. We will prove $\chil(\lbb^{\otimes N}) \geq \chil(\pi^*\lbb^{\otimes N})+ o(N^n)$. The projection formula gives $\Gamma(\widetilde{X}, \pi^*\lb^{\otimes N})= \Gamma(X, \lb^{\otimes N}\otimes \pi_*\ob_{\widetilde{X}}).$ And the natural injection $\Gamma(X, \lb^{\otimes N}) \rightarrow \Gamma(\widetilde{X}, \pi^*\lb^{\otimes N})$ is an isometry to its image under $L^2$-norms. The task is to bound the quotient.

Pick a hermitian line bundle $\mbb$ over $X$ satisfying the following two conditions:

(1) $\lbb\otimes\mbb$ is arithmetically ample, and $\lb\otimes \mb$ is ample in the classical sense;

(2) There exists an effective section $s\in \Gamma(X, \mb)$ which does not vanish at any associated point of the coherent sheaf $\pi_*\ob_{\widetilde{X}}/\ob_X$ over $X$.

The finiteness of $\pi$ implies that $\pi^*(\lbb\otimes \mbb)$ is arithmetically ample, and $\pi^*(\lb\otimes \mb)$ is ample in the classical sense. So $\chil(\pi^*(\lbb\otimes \mbb)^{\otimes N})$ and $\chil((\lbb\otimes \mbb)^{\otimes N})$ satisfy the arithmetic Hilbert-Samuel formula.

Consider the following commutative diagram of exact sequences
\[\minCDarrowwidth18pt
\begin{CD}
0 @>>> (\lb\otimes \mb)^{\otimes N} @>>> (\lb\otimes \mb)^{\otimes N} \otimes \pi_*\ob_{\widetilde{X}} @>>>
(\lb\otimes \mb)^{\otimes N}\otimes (\pi_*\ob_{\widetilde{X}}/\ob_X) @>>> 0 \\
@. @AAs^{N}A  @AAs^{N}A @AAs^{N}A @.\\
0 @>>> \lb^{\otimes N} @>>> \lb^{\otimes N}\otimes \pi_*\ob_{\widetilde{X}} @>>>
\lb^{\otimes N}\otimes (\pi_*\ob_{\widetilde{X}}/\ob_X) @>>> 0
\end{CD}\]
which will induce a diagram for long exact sequences of cohomology groups over $X$. By the choice of $s$, the three vertical morphisms are injective. Thus the diagram implies an injection
\begin{eqnarray}
\Gamma(\lb^{\otimes N}\otimes \pi_*\ob_{\widetilde{X}})/\Gamma(\lb^{\otimes N}) \rightarrow
\Gamma((\lb\otimes \mb)^{\otimes N} \otimes \pi_*\ob_{\widetilde{X}})/ \Gamma((\lb\otimes \mb)^{\otimes N})
\end{eqnarray}
which is norm-contractive. Here we endow both quotient modules with the quotient norms induced from the $L^2$-norms of global sections.

Since $\pi^*(\lb\otimes \mb)$ is ample, the section ring $\bigoplus_{N=0}^{\infty}\Gamma(\widetilde{X}, \pi^*(\lb\otimes \mb)^{\otimes N})$ is a finitely generated $\ZZ$-algebra. By picking a set of generators, one sees that
there exists a constant $c>0$ such that
$\Gamma(X, (\lb\otimes \mb)^{\otimes N} \otimes \pi_*\ob_{\widetilde{X}}) = \Gamma(\widetilde{X}, \pi^*(\lb\otimes \mb)^{\otimes N})$
is generated by sections with norms less than $c^N$. And thus $\Gamma((\lb\otimes \mb)^{\otimes N} \otimes \pi_*\ob_{\widetilde{X}})/ \Gamma((\lb\otimes \mb)^{\otimes N})$ is generated by elements with norms less than $c^N$. Applying Proposition \ref{2.1} (6) (a) to the injection (2), we have
$$
\chi\left(\Gamma(\lb^{\otimes N}\otimes \pi_*\ob_{\widetilde{X}})/\Gamma(\lb^{\otimes N})\right) \leq \chi\left(\Gamma((\lb\otimes \mb)^{\otimes N} \otimes \pi_*\ob_{\widetilde{X}})/ \Gamma((\lb\otimes \mb)^{\otimes N})\right)+O(N^{n-1}\log N).
$$
By Proposition \ref{2.1} (4) (a), we obtain
\begin{eqnarray*}
&& \chil(\pi^*\lbb^{\otimes N})-\chil(\lbb^{\otimes N})\\
&=& \chi\left(\Gamma(\lb^{\otimes N}\otimes \pi_*\ob_{\widetilde{X}})/\Gamma(\lb^{\otimes N})\right)+O(N^{n-1}\log N) \\
&\leq& \chi\left(\Gamma((\lb\otimes \mb)^{\otimes N} \otimes \pi_*\ob_{\widetilde{X}})/ \Gamma((\lb\otimes \mb)^{\otimes N})\right)+O(N^{n-1}\log N)\\
&=&\chil(\pi^*(\lbb\otimes \mbb)^{\otimes N})-\chil((\lbb\otimes \mbb)^{\otimes N})+O(N^{n-1}\log N)\\
&=&o(N^n).
\end{eqnarray*}
Here the last equality holds since $\chil(\pi^*(\lbb\otimes \mbb)^{\otimes N})$ and $\chil((\lbb\otimes \mbb)^{\otimes N})$ have the Hilbert-Samuel formula with the same leading term.
\end{proof}

\

\subsubsection{A Reduction}\label{reduction}

Keep the notation in Theorem \ref{main}. We claim that it suffices to prove the inequality under the following three assumptions:
\begin{enumerate}
\item[(1)] $X$ is normal and generically smooth.
\item[(2)] $\lbb$ and $\mbb$ are ample with positive curvatures.
\item[(3)] There is a section $s\in \gs{\mb}$ such that

(a) $s$ is effective, i.e., $\supnorm{s} \leq 1$;

(b) Each component of the Weil divisor $\mathrm{div}(s)$ is a Cartier divisor.
\end{enumerate}

We can reduce the problem to (1) by Proposition \ref{2.2}. Fix a generic resolution $\pi: \widetilde{X}\rightarrow X$, where $\widetilde{X}$ is normal and generically smooth. Replace the problem $(X, \lbb, \mbb, \ebb)$ by
$(\widetilde{X}, \pi^*\lbb,\pi^*\mbb,\pi^*\ebb)$.

Next, we see that we can replace $(\lbb, \mbb, \ebb)$ by the tensor power $(r\lbb, r\mbb, \ebb)$ for any positive integer $r$. Actually, the result for $(r\lbb, r\mbb, \ebb)$ gives
$$\chi_{\sup}(\ebb+Nr\lbb- Nr\mbb)
\geq \frac{\chern(\lbb)^n-n\cdot\chern(\lbb)^{n-1}\chern(\mbb)}{n!}(Nr)^n +o(N^n).$$
Replacing $\ebb$ by $\ebb+k\lbb-k\mbb$ for $k=0, 1, \cdots, r-1$ in the above, the desired bound for $(\lbb, \mbb)$ is obtained.

Now we can reduce the problem to (2). Assuming the result is true under (2), we need to extend the result to arbitrary ample hermitian line bundles $\lbb, \mbb$ (whose curvatures are only semipositive in general). Pick an ample hermitian line bundle $\overline{\mathcal{H}}$ over $X$ with positive curvature. Let $\epsilon$ be a positive rational number. Then $\lbb+\epsilon\overline{\mathcal{H}},\mbb+\epsilon\overline{\mathcal{H}}$ are ample with positive curvature in the sense of $\QQ$-divisors. The result for 
$(\widetilde{X}, \lbb+\epsilon\overline{\mathcal{H}},\mbb+\epsilon\overline{\mathcal{H}},\ebb)$
gives
$$\chi_{\sup}(\ebb+N\lbb- N\mbb)   \geq \frac{\chern(\lbb+\epsilon\overline{\mathcal{H}})^n-n\cdot\chern(\lbb+\epsilon\overline{\mathcal{H}})^{n-1}\chern(\mbb+\epsilon\overline{\mathcal{H}})}{n!}N^n +o(N^n).$$
We rewrite it as 
$$\chi_{\sup}(\ebb+N\lbb- N\mbb)
\geq \frac{\chern(\lbb)^n-n\cdot\chern(\lbb)^{n-1}\chern(\mbb)+\alpha(\epsilon)}{n!}N^n +o(N^n),$$
where 
$$\alpha(\epsilon)=\chern(\lbb+\epsilon\overline{\mathcal{H}})^n-n\cdot\chern(\lbb+\epsilon\overline{\mathcal{H}})^{n-1}\chern(\mbb+\epsilon\overline{\mathcal{H}}) 
- \left(\chern(\lbb)^n-n\cdot\chern(\lbb)^{n-1}\chern(\mbb) \right) $$
is a polynomial in $\epsilon$ whose constant term is 0.

We plan to obtain the expected result by taking $\epsilon \rightarrow 0$. The difficulty is that the error term $o(N^n)$ actually depends on $\epsilon$. The following lemma solves this little problem. 
\begin{lem}\label{regularity}
Let $F, G, H$ be three functions of $\mathbb{N} \rightarrow \RR$ with $H(N)>0$ for all $N\in \mathbb{N}$, and let $\alpha: \QQ_{>0} \rightarrow \RR$ be a function such that $\alpha(\epsilon)\rightarrow 0$ as $\epsilon\rightarrow 0.$ Assume that for any $\epsilon\in \QQ_{>0},$
$$F(N) \geq G(N)+\alpha(\epsilon)H(N)+o_\epsilon(H(N)),  \quad N\rightarrow \infty,$$ 
where the error term $o_\epsilon(H(N))$ depends on $\epsilon$.
Then we have
$$F(N) \geq G(N)+o(H(N)),  \quad N\rightarrow \infty.$$
\end{lem}
\begin{proof}
It is an exercise in Calculus. For an explanation of the inequalities, see the remark after Theorem \ref{main}.
By dividing the inequalities by $H(N)$, we can assume that $H(N)=1$. By replacing $F$ by $F-G$, we can assume that $G=0.$

Now we see that there exists a function $R(\epsilon, N)$ such that $\lim_{N\rightarrow \infty} R(\epsilon, N)=0$ and
$$F(N) \geq \alpha(\epsilon)+R(\epsilon, N).$$
We want to show that 
$$F(N) \geq o(1),  \quad N\rightarrow \infty,$$
which says that $F$ is greater than or equal to some function on $\mathbb{N}$ with limit 0 at infinity. 

Use proof by contradiction. Assume that it is not true. Then there exists a real number $c>0$, such that $F(N)<-c$ holds for infinitely many $N$. 

One can find an $\epsilon_0>0$ such that $\displaystyle\alpha(\epsilon_0)>-\frac{c}{2}.$ We also have 
$\displaystyle R(\epsilon_0, N)>-\frac{c}{2}$ when $N$ is large enough. Thus
$$F(N) \geq \alpha(\epsilon_0)+R(\epsilon_0, N)  > -\frac{c}{2}-\frac{c}{2}=-c$$
for $N$ large enough. It contradicts to the above statement that $F(N)<-c$ holds for infinitely many $N$. 
\end{proof}

Now we go back to the assumption in (3). Condition (b) is automatic if $X$ is regular, which can be achieved in (1) if $X$ has a resolution of singularity. In partcular, it works if $\dim X=2$ since the resolution of singularity is proved in Lipman \cite{Li}.

In the case $\dim X>2$, we apply the arithmetic Bertini theorem proved by Moriwaki \cite{Mo1}. See Page 1326 of the paper. By (1), we can assume that $X$ is generically smooth. We claim that there exists a (non-empty) open subscheme $U$ of $\mathrm{Spec}(\ZZ)$ such that $X_U$ is smooth over $U$. In fact, we can find an open subset of $X$ which is smooth and contains the generic fibre since the property of being smooth is open. Let $V$ be the complement of this open subset in $X$. Then the image of $V$ in $\mathrm{Spec}(\ZZ)$ doesn't contain the generic point, and is closed by the properness of $X$. We take $U$ to be the complement of this image in $\mathrm{Spec}(\ZZ)$. 

The closed subset $X-X_U$ 
consists of finitely many vertical fibres of $X$, and hence finitely many irreducible components. 
Denote their generic points by $\eta_1, \cdots, \eta_t$. By Moriwaki's arithmetic Bertini theorem, there exists an effective section $s\in \gs{r\mb}$ for some integer $r>0$ such that $\divv(s)_\QQ$ is smooth over $\QQ$ and $s$ does not vanish at any of $\eta_i$. We will show that $s$ satifies (b) so that (3) is achieved by replacing $(\lbb, \mbb)$ by $(r\lbb, r\mbb)$. 

Any vertical component of $\mathrm{div}(s)$, which is necessarily of codimension one, is not equal to any $\eta_i$. Therefore, it must be contained in the regular open subscheme $X_U$, and is a Cartier divisor. 

Let $K$ be the unique number field such that $X_{\QQ}\rightarrow \Spec(\QQ)$ factors through $\Spec(K)$ and $X_{\QQ}\rightarrow \Spec(K)$ is geometrically connected. By \cite[Corollary 7.9, Chapter III]{Ha}, $\divv(s)_\QQ\times_K \overline K$ is connected, and thus $\divv(s)_\QQ$ is connected. 
But $\divv(s)_\QQ$ is smooth over $\QQ$, it must be irreducible and reduced. It follows that $\mathrm{div}(s)$ has only one horizontal component with multiplicity one, which is forced to be a Cartier divisor since any other components are Cartier divisors.

\

\subsection{Analytic Parts}\label{analytic}
In this section, we prove a volume comparison theorem in the first two subsections and show a Gromov type of norm inequality in the third subsection. This section is divided into three subsections according to different settings.

Suppose $X$ is a generically smooth arithmetic variety of dimension $n$. Let $\ebb$, $\lbb$, $\mbb$, $\mbb'$ be four hermitian line bundles over $X_\CC$. We assume that 
\begin{itemize}
\item $\lbb$, $\mbb$ have positive curvatures at infinity;
\item $s\in \Gamma(\mb')$ is an effective section in the sense that $\supnorm{s}\leq 1$;
\item $X_\CC$ is endowed with the probability measure
$$\displaystyle d\mu =d\mu_{\lbb}=\frac{1}{\deg_{\lb_\CC}(X_\CC)}c_1(\lbb)^{n-1}
=\frac{1}{\int_{X_\CC} c_1(\lbb)^{n-1}}c_1(\lbb)^{n-1}$$
 induced by $\lbb$. 
\end{itemize}

Tensoring by $s$ defines an injection
$$ \Gamma(\eb+N\mathcal{L}-j\mathcal{M}) \hookrightarrow  \Gamma(\eb+N\mathcal{L}-j\mathcal{M}+\mb'). $$
\noindent This gives an induced quadratic norm
$$\left\|t\right\|_s=\left(\int |s(z)|^2|t(z)|^2d\mu\right)^{1/2}$$
for any $t\in \Gamma(\eb+N\mathcal{L}-j\mathcal{M})$.
We will compare this norm with the original $L^2$-norm
$$\left\|t\right\|_{L^2}=\left(\int |t(z)|^2d\mu\right)^{1/2}$$
on $\Gamma(\eb+N\mathcal{L}-j\mathcal{M})$. Those two norms induces $\chi_s$ and $\chi_{_{L^2}}$ on $\Gamma(\eb+N\mathcal{L}-j\mathcal{M})$. It follows from the definition that 
$$\chi_{_{L^2}}(\ebb+N\lbb-j\mbb)-\chi_s(\ebb+N\lbb-j\mbb)= \log \frac{\vol{B_{L^2}}}{\vol{B_s}}.$$
Here $B_s$ and $B_{L^2}$ are the unit balls in $\Gamma(\eb+N\mathcal{L}-j\mathcal{M})_\RR$ corresponding to these two norms, and $\vol{B_{L^2}}/\vol{B_s}$ is independent of the Haar measure chosen on $\Gamma(\eb+N\mathcal{L}-j\mathcal{M})_\RR$. The main result in this section is as follows:

\begin{thm}\label{3.1}
As $N,j \rightarrow \infty$,
$$ \chi_{_{L^2}}(\ebb+N\lbb-j\mbb)-\chi_s(\ebb+N\lbb-j\mbb) \geq  \rank \Gamma(N\mathcal{L})  \left(\int \log |s(z)| d\mu\right)   (1+O(\frac{1}{j}+\frac{1}{N})).$$
\end{thm}

\subsubsection{Distortion Functions}
In this subsection, $X$ denotes a compact complex manifold of complex dimension $n-1$ with a probability measure $d\mu$.
Here $X$ is not necessarily connected, but we require that each connected component is of dimension $n-1$. Many results over connected manifolds can be extended naturally to this case.

For any hermitian line bundle $\lbb$  over $X$, the $L^2$-norm makes $\Gamma(\lb)$ a complex Hilbert space. Suppose $s_1, s_2, \cdots, s_r $ form an orthonormal basis. Define the distortion function $b(\lbb):X\rightarrow \textbf{R}$ by
$$b(\lbb)(z)=|s_1(z)|^2+|s_2(z)|^2+\cdots+|s_r(z)|^2$$
which is independent of the basis. For convenience, define $b(\lbb)$ to be zero everywhere if $\Gamma(\mathcal{L})=0$. The following theorem was proved independently by Bouche \cite{Bo} and Tian \cite{Ti}:

\begin{thm}\label{3.2}
If $\lbb$ has positive curvature and the measure $d\mu$ over $X$ is induced by $\lbb$, then for any hermitian line bundle $\ebb$,
$$b(\ebb+N\lbb)(z)=\dim\Gamma(\eb+N\mathcal{L})(1+O(\frac{1}{N}))$$ uniformly on $X$ as $N \rightarrow \infty$.
\end{thm}

Now we generalize it to an estimate on $\ebb+N\lbb-j\mbb$ for hermitian line bundles $\ebb$, $\lbb$, $\mbb$ over $X$, which will be used to prove Theorem \ref{3.1}.

\begin{pro}\label{3.3}
If $\lbb$ and $\mbb$ have positive curvatures, and the measure $d\mu$ over $X$ is induced by $\lbb$, then for any hermitian line bundle $\ebb$,
$$b(\ebb+N\lbb-j\mbb)(z) \leq \dim\Gamma(\eb+N\mathcal{L})(1+O(\frac{1}{N}+\frac{1}{j}))$$ uniformly on $X$ as $N, j \rightarrow \infty$.
\end{pro}

\begin{proof}
Assume $\ebb$ to be trivial as usual. For fixed $z\in X$, one can choose an orthonormal basis of $\Gamma(j\mb)$ under the measure $d\mu_{\mbb}$ such that only one section in this basis is nonzero at $z$. Call this section $s^j$. Then by Theorem \ref{3.2}, one has

$|s^j(z)|^2=b(j\mbb)(z)=\dim\Gamma(j\mb)(1+O(1/j)),$

$\supnorm{s^j}^2=\sup_{x\in X}|s^j(x)|^2 \leq \sup_{x\in X} b(j\mbb)(x) = \dim\Gamma(j\mb)(1+O(1/j)).$

\noindent Those imply $|s^j(z)|/\supnorm{s^j}=1+O(1/j).$ Note that this result actually does not depend on the measure on $X$. Next we use the measure $d\mu=d\mu_{\lbb}$.

For each such $s^j$, consider the two quadratic norms $\lnorm{\cdot}$ and $\|\cdot\|_{s^j}$ on $\Gamma(N\lb-j\mb)$. By linear algebra, there exists a basis $t_1, t_2, \cdots, t_r$, which is orthonormal under $\lnorm{\cdot}$ and orthogonal under $\|\cdot\|_{s^j}$. Since $\|\cdot\|_{s^j}$ is the induced norm under
$$ \Gamma(N\mathcal{L}-j\mathcal{M}) \stackrel{s^j}{\rightarrow}  \Gamma(N\mathcal{L}), $$
we can view $s^jt_1, s^jt_2, \cdots, s^jt_r$ as $r$ orthogonal elements of $\Gamma(N\lb)$. Normalize them and apply Theorem \ref{3.2} again:
$$\sum_{k=1}^r \frac{|s^j(z)t_k(z)|^2}{\lnorm{s^jt_k}^2} \leq   b(N\lbb)(z) = \dim\Gamma(N\lb)(1+O(\frac{1}{N})).$$
Since
$$\lnorm{s^jt_k}^2=\int |s^j(x)t_k(x)|^2 d\mu \leq \supnorm{s^j}^2 \int |t_k(x)|^2 d\mu=\supnorm{s^j}^2=|s^j(z)|^2(1+O(\frac{1}{j})),$$
we have
$$\sum_{k=1}^r \frac{|s^j(z)t_k(z)|^2}{|s^j(z)|^2(1+O(\frac{1}{j}))} \leq   
 \dim\Gamma(N\lb)(1+O(\frac{1}{N})).$$
It simplifies to 
$$\sum_{k=1}^r |t_k(z)|^2 \leq   \dim\Gamma(N\lb)(1+O(\frac{1}{N}+\frac{1}{j})),$$
which is exactly 
$$b(N\lbb-j\mbb)(z) \leq \dim\Gamma(N\mathcal{L})(1+O(\frac{1}{N}+\frac{1}{j})).$$
\end{proof}

\begin{remark}
Since any hermitian line bundle is the difference of two ample hermitian line bundles with positive curvatures, this result actually gives an upper bound of $b(N\lbb)$ for any hermitian line bundle $\lbb$.
\end{remark}

\subsubsection{The Comparison}
Now we can prove Theorem \ref{3.1}. We follow the same strategy as \cite[Lemma 3.8]{AB}.
\begin{proof}[Proof of Theorem \ref{3.1}]
As in the proof of Proposition \ref{3.3}, pick an $\RR$-basis $t_1, t_2, \cdots, t_r$ of $\Gamma(\eb+N\lb-j\mb)_\RR$ which is orthonormal under $\lnorm{\cdot}$ and orthogonal under $\|\cdot\|_s$. Then
$$
\log(\frac{\vol{B_{L^2}}}{\vol{B_s}})=\log\prod_{k=1}^r \|t_k \|_s=\frac{1}{2}\sum_{k=1}^r\log \int |s(x)|^2\cdot |t_k(x)|^2 d\mu.
$$
Since $\displaystyle\int |t_k(x)|^2 d\mu=1$, one can view $|t_k(x)|^2 d\mu$ as a probability measure on $X_\CC$. Applying Jensen's inequality to the function $\log$, one gets
$$\log \int |s(x)|^2\cdot |t_k(x)|^2 d\mu \geq \int \log |s(x)|^2\cdot |t_k(x)|^2 d\mu,$$
and thus
\begin{eqnarray*}
\log(\frac{\vol{B_{L^2}}}{\vol{B_s}})
&\geq& \frac{1}{2}\sum_{k=1}^r \int \log |s(x)|^2\cdot |t_k(x)|^2 d\mu  \\
&=& \frac{1}{2}\int \log |s(x)|^2\cdot \sum_{k=1}^r |t_k(x)|^2 d\mu \\
&\geq& \frac{1}{2}\dim\Gamma(\eb+N\mathcal{L})  \left(\int \log |s(z)|^2 d\mu\right)   (1+O(\frac{1}{j}+\frac{1}{N}))\\
&=& \dim\Gamma(N\mathcal{L})  \left(\int \log |s(z)|d\mu\right)   (1+O(\frac{1}{j}+\frac{1}{N})),
\end{eqnarray*}
where the last inequality uses Proposition \ref{3.3} and the assumption $\supnorm{s}\leq 1$.
\end{proof}

\subsubsection{Gromov's Inequality}

To end this section, we show a version of Gromov's norm comparison theorem for there line bundles. The proof is the same as the original one in \cite{GS2}. We still include it here.

\begin{pro}\label{3.4}
Suppose $X$ is a compact complex manifold of complex dimension $n$ endowed with a volume form $d\mu$ which is positive everywhere. Let $\lbb$, $\mbb$ and $\ebb$ be three hermitian line bundles over $X$. Then there exists a positive constant $c$ such that
$$ \lnorm{s} \geq c(k+j)^{-n}\supnorm{s} \ ,  \ \forall k,j>0, s\in\Gamma(\eb+k\mathcal{L}+j\mathcal{M}). $$
\end{pro}

\begin{proof}
We only consider the case where $\ebb$ is the trivial hermitian line bundle. One can find a finite open cover $\{U_\alpha\}_\alpha$ of X satisfying the following

(1) $\{U_\alpha\}_\alpha$ trivializes $\lb$ and $\mb$;

(2) Each $U_\alpha$ is isomorphic to the polydisc $\{z\in \CC^n: |z|<3\}$ under the coordinate $z_\alpha: U_\alpha \rightarrow \CC^n$;

(3) The discs $\{x\in U_\alpha: |z_\alpha(x)|<1\}$ defined by these coordinates still cover $X$.

\noindent Any section $s\in\Gamma(k\mathcal{L}+j\mathcal{M})$ corresponds to a set of holomorphic functions $\{s_\alpha: U_\alpha \rightarrow \CC\}_\alpha$ under the trivialization by $\{U_\alpha\}_\alpha$. We also view $s_\alpha$ as a holomorphic function on the polydisc $\{z\in \CC^n: |z|<3\}$ by the coordinate $z_\alpha$.

Suppose $h_\alpha$, $h'_\alpha$ give the metrics of $\lbb$ and $\mbb$. More precisely, $h_\alpha$ and $h'_\alpha$ are  infinitely-differentiable functions on $U_\alpha \rightarrow \RR_{>0}$ such that the metric of $s=\{s_\alpha\}_\alpha \in\Gamma(k\mathcal{L}+j\mathcal{M})$ is given by $|s|^2=h_\alpha^k h_\alpha'^j s_\alpha \overline{s}_\alpha$ in open set $U_\alpha$.

View $D_\alpha=\{z_\alpha\in \CC^n: |z_\alpha|\leq 2\}$ as a ball in $\RR^{2n}$ and $h_\alpha$ as a function on it. We can find a constant $c$ bounding the length of the gradient of $h_\alpha$ in $\{z_\alpha\in \CC^n: |z_\alpha|\leq 2\}$ for all $\alpha$.
Pick a constant $c_1>1$ such that $\displaystyle c_1> \max\left\{\frac{c}{h_\alpha(x)},\frac{c}{h'_\alpha(x)}\right\}, \forall x \in D_\alpha, \forall\alpha$.

For any $x_0, x_1\in D_\alpha$, \ one has $|h_\alpha(x_1)- h_\alpha(x_0)|\leq c|z_\alpha(x_1)-z_\alpha(x_0)|$ and thus
$$ h_\alpha(x_1)\geq h_\alpha(x_0)-c|z_\alpha(x_1)-z_\alpha(x_0)| \geq h_\alpha(x_0)(1-c_1|z_\alpha(x_1)-z_\alpha(x_0)|)$$

Now consider the norms of $s=(s_\alpha)\in\Gamma(k\mathcal{L}+j\mathcal{M})$. Suppose $\supnorm{s}=|s(x_0)|$ for a point $x_0\in X$. Suppose $x_0$ is contained in $\{x\in U_\alpha: |z_\alpha(x)|<1\}$. Now the neighborhood $U=\{x\in U_\alpha: |z_\alpha(x)-z_\alpha(x_0)|<c_1^{-1}\}$ is contained in $D_\alpha$ since $c_1>1$. We have $h_\alpha(x)\geq h_\alpha(x_0)(1-c_1|z_\alpha(x)-z_\alpha(x_0)|)$ for all $x\in U$ and the same for $h'_\alpha$.

Finally, we come to our estimate. For simplicity, assume $z_\alpha(x_0)=0.$ Then
\begin{eqnarray*}
&& \lnorm{s}^2 \\
&\geq& \int_U |s(x)|^2 d\mu \\
&\geq& c' \int_{\{z\in\CC^n: \  |z|<c_1^{-1}\}} h_\alpha(z)^k h'_\alpha(z)^j s_\alpha(z)\overline{s_\alpha(z)} dV(z)\\
&\geq& c' \int_{\{z\in\CC^n: \  |z|<c_1^{-1}\}} h_\alpha(0)^k h'_\alpha(0)^j (1-c_1|z|)^k (1-c_1|z|)^j s_\alpha(z)\overline{s_\alpha(z)} dV(z)\\
&=& c''h_\alpha(0)^k h'_\alpha(0)^j \int_0^{c_1^{-1}} \int_{\{z\in\CC^n: \  |z|=1\}}(1-c_1|rz|)^{k+j} s_\alpha(rz)\overline{s_\alpha(rz)}  r^{2n-1}dS(z)dr \\
&=& c''h_\alpha(0)^k h'_\alpha(0)^j \int_0^{c_1^{-1}} \left(\int_{\{z\in\CC^n: \  |z|=1\}}s_\alpha(rz)\overline{s_\alpha(rz)} dS(z)\right) (1-c_1r)^{k+j} r^{2n-1}dr.
\end{eqnarray*}
Here a few things need to be explained. The measure $dV(z)$ is the standard one for the ball, and $dS(z)$ is the Lebesgue measure for the unit sphere. The constants $c'$ is such that $d\mu(z) \geq c' \ dV(z)$, which exists because $d\mu$ is pointwise positive. The constant $c''$ occurs by the spherical coordinate, which asserts that $dV(z)$ is equal to some constant multiple of $r^{2n-1}dS(z)dr$.

Now we integrate $s_\alpha(rz)\overline{s_\alpha(rz)}$ over the unit sphere. Since the function $s_\alpha(rz)\overline{s_\alpha(rz)}$ is pluri-subharmonic, the mean value inequality is valid. Therefore,
\begin{eqnarray*}
 \lnorm{s}^2 
&\geq& c''h_\alpha(0)^k h'_\alpha(0)^j \int_0^{c_1^{-1}}s_\alpha(0)\overline{s_\alpha(0)}(1-c_1r)^{k+j} r^{2n-1}dr\\
&=& c''c_1^{-2n}|s(x_0)|^2 \int_0^1 (1-r)^{k+j} r^{2n-1}dr\\
&=& c''c_1^{-2n}\frac{\supnorm{s}^2}{(2n+k+j){2n+k+j-1 \choose 2n-1}},
\end{eqnarray*}
which implies our result.
\end{proof}

\subsection{Proof of the Bigness Theorems}\label{proof}

The task of this section is to prove Theorem \ref{main} in the smooth case proposed by Section \ref{reduction}. 
The proof we are giving here is analogous to the ample case in \cite{AB}.
As a byproduct, we prove Theorem \ref{1.3}.

\subsubsection{Proof of Theorem \ref{main}}
Before the proof of Theorem \ref{main}, we state a key result in a slightly different setting. Let $\ebb$, $\lbb$, $\mbb$, $\mbb'$ be four hermitian line bundles over an arithmetic variety $X$ of dimension $n$. Assume that:
\begin{itemize}
\item $X$ is normal and generically smooth;
\item $\lbb$, $\mbb$ are ample with positive curvatures;
\item There is a nonzero section $s\in \Gamma(\mbb')$ such that $\supnorm{s} < 1$ and each component of the Weil divisor $\mathrm{div}(s)$ is a Cartier divisor. 
\end{itemize} 
\begin{pro}\label{4.1}
As $N,j\rightarrow \infty$,
$$\chil(\ebb+N\lbb- j\mbb)-\chil(\ebb+N\lbb-j\mbb+\mbb')
\geq  -\frac{\chern(\lbb)^{n-1}\chern(\mbb')}{(n-1)!}N^{n-1}+ O(N^{n-2}(\frac Nj+\log j)).
$$
\end{pro}

\

The situation is similar to that in Section \ref{analytic}. Recall that we have the injection
$$ \Gamma(\eb+N\mathcal{L}-j\mathcal{M}) \stackrel{\otimes s}{\hookrightarrow } \Gamma(\eb+N\mathcal{L}-j\mathcal{M}+\mb').$$
And we also have two quadratic norms 
$$\left\|t\right\|_{L^2}=\left(\int |t(z)|^2d\mu\right)^{1/2}$$
and
$$\left\|t\right\|_s=\left(\int |s(z)|^2|t(z)|^2d\mu\right)^{1/2}$$
on $\Gamma(\eb+N\mathcal{L}-j\mathcal{M})$. Here $d\mu$ is the probability measure induced by $\lbb$.
They define two arithmetic volumes $\chil$ and  $\chi_s$ on $\Gamma(\eb+N\mathcal{L}-j\mathcal{M})$. Then Theorem \ref{3.1} asserts
$$ \chi_{_{L^2}}(\ebb+N\lbb-j\mbb)-\chi_s(\ebb+N\lbb-j\mbb) \geq  \mathrm{rank} \Gamma(N\mathcal{L})  \left(\int \log |s(z)| d\mu\right) +O(\frac{N^{n-1}}{j}+N^{n-2}).$$

Another estimate needed for Proposition \ref{4.1} is:
\begin{lem}\label{4.3}
If furthermore $\divv(s)$ is a prime divisor, then
$$\chi_s(\ebb+N\lbb-j\mbb)-\chil(\ebb+N\lbb-j\mbb+\mbb') \geq - \frac{\chern(\lbb)^{n-1}\cdot \divv(s)}{(n-1)!}N^{n-1}+O(N^{n-2}\log (N+j)).$$
\end{lem}

We will prove this lemma later, but we will first consider its consequences.

\

\begin{proof}[Proof of Proposition \ref{4.1}]
\

Proposition \ref{4.1} is a summation of Theorem \ref{3.1} and Lemma \ref{4.3}. We first consider the case that $\divv(s)$ is a prime divisor, i.e., it satisfies the extra condition of Lemma \ref{4.3}. Then
\begin{eqnarray*}
&& \chil(\ebb+N\lbb- j\mbb)-\chil(\ebb+N\lbb-j\mbb+\mbb') \\
&\geq&  -\frac{\chern(\lbb)^{n-1}\cdot \divv(s)}{(n-1)!}N^{n-1}+ \mathrm{rank} \Gamma(N\mathcal{L})  \left(\int \log |s(z)| d\mu\right)+O(\frac{N^{n-1}}{j}+N^{n-2}\log (N+j))\\
&=& -\frac{\chern(\lbb)^{n-1}\cdot \divv(s)}{(n-1)!}N^{n-1}+ \frac{N^{n-1}\deg_{\lb_{\QQ}}(X_\QQ)}{(n-1)!}
\int \log |s(z)| d\mu +O(N^{n-2}(\frac Nj+\log j))\\
&=& -\frac{\chern(\lbb)^{n-1}\cdot \divv(s)}{(n-1)!}N^{n-1}+ \frac{N^{n-1}}{(n-1)!}
\int \log |s(z)| c_1(\lbb)^{n-1} +O(N^{n-2}(\frac Nj+\log j))\\
&=& -\frac{N^{n-1}}{(n-1)!} \left(\chern(\lbb)^{n-1}\cdot \divv(s)-
\int \log |s(z)| c_1(\lbb)^{n-1}\right) +O(N^{n-2}(\frac Nj+\log j))\\
&=& -\frac{N^{n-1}}{(n-1)!} \chern(\lbb)^{n-1}\cdot (\divv(s), -\log |s|^2) +O(N^{n-2}(\frac Nj+\log j))\\
&=& -\frac{N^{n-1}}{(n-1)!} \chern(\lbb)^{n-1}\chern(\mbb') +O(N^{n-2}(\frac Nj+\log j)).
\end{eqnarray*}
Here the error term 
$$\frac{N^{n-1}}{j}+N^{n-2}\log (N+j)
=N^{n-2}(\frac Nj+\log (1+\frac Nj)+\log j)=O(N^{n-2}(\frac Nj+\log j)).$$

Now we consider the general case.  By the assumption, there exists $s\in \Gamma(\mb')$ with $\supnorm{s} < 1$, and we can decompose $\mb'=\mb_1+\cdots+\mb_r$ such that $s=s_1\otimes \cdots \otimes s_r$ for $s_k\in \Gamma(\mb_k)$ and each $\divv(s_k)$ is a prime divisor. We claim that we can endow each $\mb_k$ with a hermitian metric $\|\cdot\|_k$ such that 
$\|s\|=\prod_{k=1}^{r} \|s_k\|_k$ and each $\|s_k\|_{k,\sup} < 1$. Once this is true, we have 
\begin{eqnarray*}
&& \chil(\ebb+N\lbb- j\mbb+\mbb_1+\cdots+\mbb_{k-1})-\chil(\ebb+N\lbb-j\mbb+\mbb_1+\cdots+\mbb_k)\\
&\geq& -\frac{N^{n-1}}{(n-1)!} \chern(\lbb)^{n-1}\chern(\mbb_{k}) +O(N^{n-2}(\frac Nj+\log j)).
\end{eqnarray*}
Take the summation for $k=1,\cdots, r$, we obtain the result.

In the end, we check the existence of the metric $\|\cdot\|_k$. Take any metric $\|\cdot\|_k$ for $\mb_k$ such that $\|s\|=\prod_{k=1}^{r} \|s_k\|_k$. Take a positive number $\epsilon>0$, define 
$\|\cdot\|'_k=\|\cdot\|_k/(\epsilon +\| s_k\|_k)$ for $k=1,2,\cdots, r-1$, and
$\|\cdot\|'_r=\displaystyle \left(\prod_{k=1}^{r-1}(\epsilon +\| s_k\|_k)\right) \ \|\cdot\|_r$. 
Then 
$$\|s_k\|'_k=\frac{\|s_k\|_k}{\epsilon +\| s_k\|_k} \leq \frac{\|s_k\|_{k,\sup}}{\epsilon +\| s_k\|_{k,\sup}}<1, \quad k=1,2,\cdots, r-1.$$
By
$\displaystyle \|s_r\|'_r =\left(\prod_{k=1}^{r-1}(\epsilon +\| s_k\|_k)\right) \ \|s_r\|_r 
=  \|s\| +O(\epsilon)$,
we see that $\|s_r\|'_{r,\sup}<1$ for $\epsilon$ small enough.
Then $\|\cdot\|'_k$ satisfies the conditions.
\end{proof}

\

\begin{proof}[Proof of Theorem \ref{main}]
Theorem \ref{main} is implied by Proposition \ref{4.1} by simple computation. Keep the notation in Theorem \ref{main}. By Section \ref{reduction}, it suffices to prove the inequality under the following three assumptions:
\begin{enumerate}
\item[(1)] $X$ is normal and generically smooth.
\item[(2)] $\lbb$ and $\mbb$ are ample with positive curvatures.
\item[(3)] There is a section $s\in \gs{\mb}$ such that

(a) $s$ is effective, i.e., $\supnorm{s} \leq 1$;

(b) Each component of the Weil divisor $\mathrm{div}(s)$ is a Cartier divisor.
\end{enumerate}

Recall that we need to prove
$$\chisup(\ebb+N(\lbb-\mbb))
\geq \frac{\chern(\lbb)^n-n\cdot\chern(\lbb)^{n-1}\chern(\mbb)}{n!}N^n +o(N^n).$$
We only need to show the above result for $\chil$ by Proposition \ref{result of gromov} (2).

Apply Proposition \ref{4.1} to the quadruple $(\ebb, \lbb, \mbb, \mbb'=\mbb)$. We have
\begin{eqnarray*}
&& \chil(\ebb+N\lbb- N\mbb)\\
&=& \chil(\ebb+N\lbb)+ \sum_{j=1}^N \left(\chil(\ebb+N\lbb- j\mbb)-\chil(\ebb+N\lbb-(j-1)\mbb)\right) \\
&\geq& \frac{\chern(\lbb)^n}{n!}N^n +O(N^{n-1}\log N) + \sum_{j=1}^N \left(-\frac{\chern(\lbb)^{n-1}\chern(\mbb)}{(n-1)!}N^{n-1}+ O(N^{n-2}(\frac Nj+\log j))\right) \\
&=& \frac{\chern(\lbb)^n-n\cdot\chern(\lbb)^{n-1}\chern(\mbb)}{n!}N^n +O(N^{n-1}\log N).
\end{eqnarray*}
Here we used the arithmetic Hilbert-Samuel formula for $\chil(\ebb+N\lbb)$. It proves the theorem. 
\end{proof}

\

\noindent It remains to prove Lemma \ref{4.3}. 

\begin{proof}[Proof of Lemma \ref{4.3}]

For simplicity of notation, we only consider the case that $\ebb$ is trivial. We need to show that 
$$\chi_s(N\lbb-j\mbb)-\chil(N\lbb-j\mbb+\mbb') \geq - \frac{\chern(\lbb)^{n-1}\cdot Y}{(n-1)!}N^{n-1}+O(N^{n-2}\log(N+j)).$$
Here we denote $Y=\divv(s)$.
The key is to analyze volume relations by the exact sequence
$$0\rightarrow \gs{N\lb-j\mb} \stackrel{\otimes s}{\rightarrow} \gs{N\lb-j\mb+\mb'} \rightarrow \Gamma(Y, N\lb-j\mb+\mb').$$

Denote $\Gamma=\Gamma(X,N\lb-j\mb+\mb')/s\Gamma(X,N\lb-j\mb)$. Then we have two exact sequences:
\begin{eqnarray*}
0&\rightarrow &\gs{N\lb-j\mb} \stackrel{\otimes s}{\rightarrow} \gs{N\lb-j\mb+\mb'} \rightarrow \Gamma \rightarrow 0,\\
0&\rightarrow &\Gamma \rightarrow \Gamma(Y, N\lb-j\mb+\mb').
\end{eqnarray*}

\noindent Two norms are induced on $\Gamma$: the quotient norm $\|\cdot\|_{q}$ and the subspace norm  $\|\cdot\|_{\mathrm{sub}}$. Let $\chi_q(\Gamma)$ and $\chi_{\mathrm{sub}}(\Gamma)$ be the corresponding arithmetic volumes.

By Proposition \ref{2.1} (4) (a), one has
$$\chi_s(N\lbb-j\mbb)-\chil(N\lbb-j\mbb+\mbb')+ \chi_q(\Gamma)\geq 0.$$
So it suffices to show
$$\chi_q(\Gamma)\leq \frac{\chern(\lbb)^{n-1}\cdot Y}{(n-1)!}N^{n-1}+O(N^{n-2}\log(N+j)).$$

We first consider the case that $Y$ is a vertical divisor, i.e., it is a component of the fibre of 
$X$ over some prime $p$ of $\ZZ$. Then $\Gamma(Y, N\lb-j\mb+\mb')$ and $\Gamma$ are torsion. It follows that
\begin{eqnarray*}
\chi_q(\Gamma) 
&=& \log\# \Gamma \leq  \log\#\Gamma(Y, N\lb-j\mb+\mb') 
\leq  \log\#\Gamma(Y, N\lb+\mb')\\
&=& \dim_{\mathbb{F}_{p}} \Gamma(Y, N\lb+\mb') \log p
= \frac{\chern(\lbb)^{n-1}\cdot Y}{(n-1)!}N^{n-1}+O(N^{n-2}).
\end{eqnarray*}

Here  $\chern(\lbb)^{n-1}\cdot Y=c_1(\lb|_{Y})^{n-1} \log p$ by definition of the arithmetic intersection, and the last equality is just the classical Hilbert-Samuel formula over the variety $Y$. 

Now we consider the case that $Y$ is a horizontal divisor. Then $Y$ is an arithmetic variety. By Lemma \ref{4.4} below, we have  
$\chi_q(\Gamma) \leq \chi_{\mathrm{sub}}(\Gamma) + O(N^{n-2}\log(N+j)). $
Thus it suffices to show that 
$$\chi_{\mathrm{sub}}(\Gamma) \leq \frac{\chern(\lbb)^{n-1}\cdot Y}{(n-1)!}N^{n-1}+O(N^{n-2}\log(N+j)).$$

We have an injection $\Gamma(Y, N\lb-j\mb+\mb') \hookrightarrow \Gamma(Y, N\lb+\mb')$ given by tensoring by any section $s'\in \Gamma(Y, j\mb)$. When $N$ is large enough, we can have $\supnorm{s'}<1$ so that the injection is norm-contractive.  By the ampleness theorem of Zhang, $\Gamma(Y, N\lb+\mb')$ is generated by elements with $L^2$-norms less than 1.
Therefore, apply Proposition \ref{2.1} (6) (a) to the injection $\Gamma \hookrightarrow \Gamma(Y, N\lb+\mb')$ and obtain
\begin{eqnarray*}
\chi_{\mathrm{sub}}(\Gamma)
&\leq&  \chil(\Gamma(Y, N\lb+\mb'))+\log V(\rank\Gamma)-\log V(\rank\Gamma(Y, N\lb+\mb'))\\
&=&  \chil(\Gamma(Y, N\lb+\mb'))+O(\rank\Gamma(Y, N\lb+\mb')\log\rank\Gamma(Y, N\lb+\mb'))\\
&=&\frac{\chern(\lbb)^{n-1}\cdot Y}{(n-1)!}N^{n-1}+O(N^{n-2}\log N),
\end{eqnarray*}
where the last equality follows from the arithmetic Hilbert-Samuel formula over $Y$.
It finishes the proof.
\end{proof}

\

\noindent The last part of the proof is the following lemma.
\begin{lem}\label{4.4}
$\chi_q(\Gamma)  \leq \chi_{\mathrm{sub}}(\Gamma)+ O(N^{n-2}\log(N+j)). $
\end{lem}

\begin{proof}
Denote the quotient map by $\phi:\Gamma(X,N\lb-j\mb+\mb')\rightarrow \Gamma$. Applying Proposition \ref{3.4}, we get for any $\gamma \in \Gamma$,
$$\|\gamma\|_{q}
= \inf_{t\in \phi^{-1}(\gamma)} \lnorm{t}
\geq c(N+j)^{-n}\inf_{t\in \phi^{-1}(\gamma)} \supnorm{t}
\geq c(N+j)^{-n}\inf_{t\in \phi^{-1}(\gamma)} \lnorm{t|_{Y}}
= c(N+j)^{-n} \|\gamma\|_{\mathrm{sub}}.$$
We have $B_q(\Gamma) \subseteq c^{-1}(N+j)^n B_{\mathrm{sub}}(\Gamma)$, and $\vol{B_q(\Gamma)} \leq (c^{-1}(N+j)^n)^{\rank\Gamma}\vol{B_{\mathrm{sub}}(\Gamma)}.$
Therefore,
$$\chi_q(\Gamma)- \chi_{\mathrm{sub}}(\Gamma)
= \log \frac{\vol{B_q(\Gamma)}}{\vol{B_{\mathrm{sub}}(\Gamma)}}
\leq(\rank\Gamma) \log(c^{-1}(N+j)^n)
=O(N^{n-2}\log (N+j)).
$$
\end{proof}

\subsubsection{Proof of Theorem \ref{1.3}}
We first show an estimate for $h^0=h^0_{\sup}$. Let $\ebb, \lbb, \mbb, \mbb'$ be as in the setting right before 
Proposition \ref{4.1}, we have:
\begin{lem}\label{4.1'}
As $N,j\rightarrow \infty$,
$$h^0(\ebb+N\lbb- j\mbb)-h^0(\ebb+N\lbb-j\mbb+\mbb')
\geq  O(N^{n-1}\log (N+j)).$$
\end{lem}

\begin{proof}
We will check that we can modify the proof of Proposition \ref{4.1} to obtain the above result in terms of $h^0$. We only need to show the modified versions of Theorem \ref{3.1} and Lemma \ref{4.3}, since Proposition \ref{4.1} is just a summation of them. We sketch the idea here.

For simplicity of notation, we only consider the case that $\ebb$ is trivial.
Recall that Theorem \ref{3.1} asserts that
$$ \chi_{_{L^2}}(N\lbb-j\mbb)-\chi_s(N\lbb-j\mbb) \geq  \rank \Gamma(N\mathcal{L})  \left(\int \log |s(z)| d\mu\right)   (1+O(\frac{1}{j}+\frac{1}{N})).$$
These two norms satisfies $\lnorm{\cdot} > \|\cdot\|_s$. 
By Proposition \ref{2.1}  (2), we have 
$$ h^0_{L^2}(N\lbb-j\mbb)-h^0_s(N\lbb-j\mbb) \geq  \chi_{_{L^2}}(N\lbb-j\mbb)-\chi_s(N\lbb-j\mbb)+O(N^{n-1}\log N).$$
Hence,
\begin{eqnarray}
h^0_{L^2}(N\lbb-j\mbb)-h^0_s(N\lbb-j\mbb) \geq  O(N^{n-1}\log N).
\end{eqnarray}

To obtain the counterpart for Lemma \ref{4.3}, we go to its proof. We first look at the proof of Lemma \ref{4.4}. We still use $\|\cdot\|_{q} \geq c(N+j)^{-n} \|\cdot\|_{\mathrm{sub}}$. Applying Proposition \ref{2.1} (3) (b) and the second inequality in Proposition \ref{2.1} (2), we get 
$$h^0_q(\Gamma) \leq h^0_{\mathrm{sub}}(\Gamma)+O(N^{n-2}\log (N+j)). $$
This is the counterpart of Lemma \ref{4.4}.

As for the proof of Lemma \ref{4.3}, to change $\chi$ to $h^0$ in every inequality, the substitution for 
Proposition \ref{2.1} (4) (a), (6) (a) is exactly Proposition \ref{2.1} (4) (b), (6) (b), and the substitution for the arithmetic Hilbert-Samuel formula is the expansion for $h^0$ in Corollary \ref{result of gromov} (1). Eventually, we obtain
\begin{eqnarray}
h^0_s(N\lbb-j\mbb)-h^0_{L^2}(N\lbb-(j-1)\mbb) \geq O(N^{n-1}\log (N+j)).
\end{eqnarray}

The sum of (3) and (4) is 
$$ h^0_{L^2}(N\lbb-j\mbb)-h^0_{L^2}(N\lbb-(j-1)\mbb) \geq  O(N^{n-1}\log (N+j)).$$
By Corollary \ref{result of gromov} (2), it gives 
$$ h^0(N\lbb-j\mbb)-h^0(N\lbb-(j-1)\mbb) \geq  O(N^{n-1}\log (N+j)).$$
This finishes the proof.
\end{proof}

\

\noindent
Now we prove Theorem \ref{1.3}.
\begin{proof}[Proof of Theorem \ref{1.3}]

The ``if" part is easy. Suppose $r \lbb=\mbb+\tbb$ with $r$ positive, $\mbb$ ample and $\tbb$ effective. We need to show that $\lbb$ is big. Let $\ebb$ be any line bundle. Pick an effective section $s\in \gs{\tb}$, the injection
$$\gs{\eb+N\mb}\rightarrow \gs{\eb+Nr\lb}$$
defined by tensoring by $s^{\otimes N}$ is norm-contractive. It follows that
$$h^0(\ebb+Nr\lbb)\geq h^0(\ebb+N\mbb)\geq \frac{\chern(\mbb)^n}{n!}N^n+o(N^n).$$
Set $\ebb=k\lbb$ for $k=0,1, \cdots, r-1$, we get
$$h^0(N\lbb)\geq \frac{1}{r^n} \frac{\chern(\mbb)^n}{ n!}N^n+o(N^n).$$
So $\lbb$ is big.

Now we show the other direction. Suppose $\lbb$ is big, and we need to show that it is the sum of an ample hermitian line bundle and an effective line bundle. 
Write $\lbb=\lbb'-\mbb$ for two ample hermitian line bundles $\lbb'$ and $\mbb$. By a similar process as in Proposition \ref{2.2}, we can assume that $(X, \lbb',\mbb)$ satisfies the three assumptions stated at the beginning of this section. 

Apply Lemma \ref{4.1'} in the case that $\ebb$ is trivial, $\mbb'=\mbb$ and $j=N+1$. We get 
$$ h^0(N\lbb'- (N+1)\mbb)-h^0(N\lbb'-(N+1)\mbb+\mbb) \geq  O(N^{n-1}\log N).$$
In terms of $\lbb$, it is just
$$ h^0(N\lbb-\mbb)-h^0(N\lbb) \geq  O(N^{n-1}\log N).$$
It follows that
$$ h^0(N\lbb-\mbb) \geq h^0(N\lbb) +  O(N^{n-1}\log N)>0$$
when $N$ is large enough. So $N\lbb-\mbb$ is effective and $N\lbb=\mbb+(N\lbb-\mbb)$ gives a desired decomposition.
\end{proof}

\section{Equidistribution Theory}
As an application of Theorem \ref{main}, some equidistribution theorems are generalized in this section. The equidistribution theory we are going to consider originated in the paper \cite{SUZ} of Szpiro-Ullmo-Zhang. They proved an equidistribution theorem \cite[Theorem 3.1]{BV} over complex analytic spaces for line bundles of positive curvatures over generically smooth arithmetic varieties, and it was extended to certain cases by Ullmo \cite[Theorem 2.4]{Ul} and Zhang \cite[Theorem 2.1]{Zh3} to prove the Bogomolov conjecture.

Recently, Chambert-Loir \cite{Ch2} defined the canonical measures over
Berkovich spaces, and proved an equidistribution theorem over the
Berkovich spaces \cite[Theorem 3.1]{Ch2}. It is a non-archimedean
analogue of Szpiro-Ullmo-Zhang's theory.

All the above results assume the strict positivity of the metrized
line bundle at the place where equidistribution is considered,
except for the case of curves in \cite{Ch2} which makes use of
Autissier's theorem. See Remark (3) of Definition \ref{1.1} for
Autissier's expansion. As we have seen in the introduction, we can
remove the strict positivity condition with the asymptotic result
in Theorem \ref{main}. We will put the two generalized results in Theorem \ref{5.1} as conclusions at different places. We also have Theorem \ref{5.2},
an algebraic version of Theorem \ref{5.1}. Our proof follows the
original idea of Szpiro-Ullmo-Zhang.

This section consists of five subsections. We state the main equidistribution theorems (Theorem \ref{5.1} and Theorem \ref{5.2}) in the first subsection, and prove them in the second (resp. third) subsection in the archimedean (resp. non-archimedean) case. In the fourth subsection, we extend Theorem \ref{5.1} to equidistribution of small subvarieties as what Baker-Ih \cite{BI} and Autissier \cite{Au2} did for the equidistribution of Szpiro-Ullmo-Zhang. In the fifth subsection, we consider the consequences of these theorems in algebraic dynamical systems. 

\subsection{A Generic Equidistribution Theorem}\label{generic}
Let $K$ be a number field, and $X$ be a projective variety over $K$. For each place $v$, denote by $K_v$ the $v$-adic completion of $K$, and by $\CC_v$ the completion of the algebraic closure $\overline K_v$ of $K_v$. We endow $K_v$ with the normalized absolute value $|\cdot|_v$, and $\CC_v$ the unique extension of that absolute value. There are two canonical analytic $v$-spaces:
\begin{enumerate}
\item The $\CC_v$-\textit{analytic space} $\xcv$ associated to the variety $X_{\CC_v}$. Namely, $\xcv$ is the usual complex analytic space $X_v(\CC)$ if $v$ is archimedean, and the Berkovich space associated to $X_{\CC_v}$ if $v$ is non-archimedean. See \cite{Be} for an introduction of Berkovich spaces. See also Zhang's simple description in Section 5.3.
\item The $K_v$-\textit{analytic space} $\xkv$ associated to the variety $X_{K_v}$. Namely, $\xkv$ is the usual complex analytic space $X_v(\CC)$ if $v$ is complex archimedean, the quotient of the usual complex analytic space $X_v(\CC)$ by the complex conjugate if $v$ is real archimedean, and the Berkovich space associated to $X_{K_v}$ if $v$ is non-archimedean.
\end{enumerate}

Both spaces are Hausdorff, compact, and finite disjoint unions of path-connected components. They are related by $\xkv=\xcv/\mathrm{Gal}(\overline K_v/K_v)$ as topological spaces.

We will state an equidistribution theorem over each of $\xcv$ and $\xkv$. We simply call the former the geometric case and the latter the algebraic case. One will see at the end of this subsection that the geometric case implies the algebraic case and that the algebraic cases over all finite extensions of $K$ imply the geometric case.

\subsubsection*{Geometric Case}
Let $K$ be a number field and $X$ be a projective variety of dimension $n-1$ over $K$.
Fix an embedding $\overline K\rightarrow \CC_v$ for each place $v$. We will consider equidistribution of small algebraic points over $\xcv$ for each place $v$.

We use the language of \textit{adelically metrized line bundles} by Zhang (\cite{Zh1}, \cite{Zh2}). Recall that an \textit{adelic metric} over a line bundle $\lb$ of $X$ is a $\CC_v$-norm $\|\cdot\|_v$ over the fibre $\lb_{\CC_v}(x)$ of each algebraic point $x\in X(\overline K)$ satisfying certain continuity and coherence conditions for each place $v$ of $K$.

The metric is \textit{semipositive} if it is the uniform limit of a sequence of metrics induced by integral models $(\mathcal{X}_j, \widetilde{\lb}_j)$ of $(X, \lb^{e_j})$ such that each $\widetilde{\lb}_j$ is a relatively semipositive arithmetic line bundle over $\mathcal{X}_j$. A metrized line bundle is \textit{integrable} if it is isometric to the difference of two semipositive metrized line bundles. The intersection number of integrable line bundles is uniquely defined by that limit process.

Fix an integrable line bundle $\lbb$ over $X$. The \textit{height} of $X$ is defined to be
$$\height(X)=\frac{\chern(\lbb)^n}{n\deg_{\lb}(X)}.$$
The \textit{height} of an algebraic point $x\in X(\overline{K})$ is defined to be
$$\displaystyle\height(x)=\frac{\chern(\lbb|_{\bar x})}{\deg(x)},$$
where $\bar x$ is the closure of $x$ in $X$, and $\deg(x)$ is the degree of the residue field of $\bar x$ over $K$.

Denote by $O(x)$ the \textit{Galois orbit} of $x$, the orbit of $x$ under the action of
$\mathrm{Gal}(\overline K /K)$. Then $O(x)$ is a set of algebraic points of order $\deg(x)$. One has
$$\displaystyle\height(x)= \frac{1}{\deg(x)}\displaystyle \sum_v\sum_{z\in O(x)}(-\log \|s(z)\|_v)$$
for any section $s\in \Gamma(X,\lb)$ which does not vanish at $\bar x$.

We can also view $O(x)$ as a finite subset of $\xcv$ for any place $v$. Define \textit{the probability measure associated to $x$} by
$$\displaystyle\mu_{v,x}=\frac{1}{\deg(x)}\sum_{z\in O(x)}\delta_z,$$
where $\delta_z$ is the Dirac measure of $z$ in $\xcv$.

Associated to $\lbb$, there is a \textit{v-adic canonical measure}
$c_1(\lbb)_v^{n-1}$ of total volume $\deg_{\lb}(X)$ over the space
$\xcv$ for any place $v$. When $v$ is archimedean, the measure
$c_1(\lbb)_v^{n-1}$ is simply the usual differential form in the
smooth case and extended to the general case by resolution of
singularities and some limit process. 
For limits of volume forms, we refer to \cite{BT}, \cite{De} in the analytic setting, 
and \cite{Ma}, \cite{Ch1}, \cite{Zh5} in the arithmetic setting. When $v$ is
non-archimedean, the canonical measure $c_1(\lbb)_v^{n-1}$ is
defined by Chambert-Loir in \cite{Ch2}. We will describe it in more
details when we prove equidistribution at non-archimedean places.

Now we recall some related definitions of equidistribution, which was stated in the introduction in dynamical case. The only difference is that the height of $X$ is not zero anymore.
\begin{enumerate}
\item A sequence $\{x_m\}_{m\geq 1}$ of algebraic points in $X(\overline K)$ is \textit{small} if $\height(x_m)\rightarrow \height(X)$ as $m\rightarrow \infty.$
\item A sequence $\{x_m\}_{m\geq 1}$ of algebraic points in $X(\overline K)$ is \textit{generic} if no infinite subsequence of $\{x_m\}$ is contained in a proper closed subvariety of $X$.
\item Let $\{x_m\}_{m\geq 1}$ be a sequence of algebraic points in $X(\overline K)$ and $d\mu$ a probability measure over the analytic space $\xcv$ for a place $v$ of $K$. We say the Galois orbits of $\{x_m\}$ are           \textit{equidistributed} with respect to $d\mu$ if the probability measure $\{\mu_{v,x_m}\}$ associated to the sequence converges weakly to $d\mu$ over $\xcv$; i.e.,
$$\frac{1}{\#O(x_m)}\sum_{x\in O(x_m)}f(x)  \rightarrow \int_{\xcv}f(x)d\mu$$
for any continuous function $f:\xcv \rightarrow \CC.$
\end{enumerate}

The equidistribution theorem in this subsection is the following:

\begin{thm}[Equidistribution of Small Points]\label{5.1}
Suppose $X$ is a projective variety of dimension $n-1$ over a number field $K$, and $\lbb$ is a metrized line bundle over $X$ such that $\lb$ is ample and the metric is semipositive. Let $\{x_m\}$ be an infinite sequence of
algebraic points in $X(\overline K)$ which is generic and small. Then for any place $v$ of $K$, the Galois orbits of the sequence $\{x_m\}$ are equidistributed in the analytic space $\xcv$ with respect to the probability measure $d\mu_v=c_1(\lbb)_v^{n-1}/\deg_{\lb}(X)$.
\end{thm}

\subsubsection*{Algebraic Case}
As in the geometric case, let $K$ be a number field, $X$ be a projective variety of dimension $n-1$ over $K$, and $\lbb$ be an integrable line bundle over $X$. We are going to consider equidistribution of small closed points over $\xkv$ for any place $v$.

View $X$ (resp. $X_{K_v}$) as a scheme of finite type over $K$ (resp. $K_v$). When we talk about points in $X$ or $X_{K_v}$ here, we always mean closed points in the corresponding schemes. Note that in the geometric case points are algebraic points. When $K_v\cong\CC$, there is no difference between closed points and algebraic points in $X_{K_v}$.

The \textit{height} of $X$ is still $$\height(X)=\frac{\chern(\lbb)^n}{n\deg_{\lb}(X)}.$$
The \textit{height} of a closed point $x\in X$ is defined to be
$$\displaystyle\height(x)=\frac{\chern(\lbb|_{x})}{\deg(x)}$$
where $\deg(x)$ is still the degree of the residue field of $x$ over $K$.

For any closed point $x\in X$, the base change $x_{K_v}$ splits into finitely many closed points in the scheme $X_{K_v}$. They form a set $O_v(x)$, called the \textit{Galois orbit} of $x$. We can also view $O_v(x)$ as a finite subset of $\xkv$. Define \textit{the probability measure associated to $x$} by
$$\displaystyle\mu_{v,x}=\frac{1}{\deg(x)}\sum_{z\in O_v(x)}\deg(z)\delta_z$$
where $\delta_z$ is the Dirac measure of $z$ in $\xcv$, and $\deg(z)$ is the degree of the residue field of
$z$ over $K_v$.

There is still a \textit{v-adic canonical measure} $c_1(\lbb)_v^{n-1}$ of total volume $\deg_{\lb}(X)$ over the space $\xkv$ for any place $v$. Actually the $v$-adic canonical measure here is just the push-forward measure of the one in the geometric case under the natural map $\xcv\rightarrow\xkv.$

With analogous definitions of small sequences, generic sequences, and equidistribution, we have the following algebraic version of Theorem \ref{5.1}:

\begin{thm}[Equidistribution of Small Points]\label{5.2}
Suppose $X$ is a projective variety of dimension $n-1$ over a number field $K$, and $\lbb$ is a metrized line bundle over $X$ such that $\lb$ is ample and the metric is semipositive. Let $\{x_m\}$ be an infinite sequence of
closed points in $X$ which is generic and small. Then for any place $v$ of $K$, the Galois orbits of the sequence $\{x_m\}$ are equidistributed in the analytic space $\xkv$ with respect to the canonical measure $d\mu_v=c_1(\lbb)_v^{n-1}/\deg_{\lb}(X)$.
\end{thm}

\subsubsection*{Equivalence}
Via the projection $\xcv\rightarrow\xkv,$ the push-forward measures of $\mu_{v,x_m}$ and $d\mu_v=c_1(\lbb)_v^{n-1}$ over $\xcv$ give exactly their counterparts over $\xkv$. Thus it is easy to see that Theorem \ref{5.1} implies Theorem \ref{5.2}.

Conversely, Theorem \ref{5.2} implies Theorem \ref{5.1}. The results of Theorem \ref{5.2} for all finite extensions $K'$ of $K$ imply the equidistribution of Theorem \ref{5.1}. In fact, considering the base change $X_{K'}$ of $X$,
Theorem \ref{5.2} implies
$$\int_{\xcv}f\mu_{v,x_m}\rightarrow\int_{\xcv}fd\mu_v,$$
for any continuous function $f:\xcv\rightarrow \CC$ that is the pull-back via $\xcv\rightarrow X_{K'_v}^{\mathrm{an}}$ of a continuous function over $X_{K'_v}^{\mathrm{an}}$ for some extension of the valuation $v$ to $K'$. Here all $X_{K'_v}^{\mathrm{an}}$ form a projective system of analytic spaces with limit $\xcv$.

A classical result says that any finite extension of $K_v$ is isomorphic to some $K'_v$ above. For a proof see, for example, Exercise 2 in Page 30 of \cite{Se}. Now it suffices to show that the vector space of all such $f$ is dense in the ring of continuous functions of $\xcv$. We need the following Stone-Weierstrass Theorem.

\begin{thmm}[Stone-Weierstrass]
Let $X$ be a compact Hausdorff space, $C(X)$ be the ring of real-valued continuous functions of $X$, and $V\subset C(X)$ be an $\RR$-vector space. Then $V$ is dense in $C(X)$ under the supremum norm if the following two conditions hold:
\begin{enumerate}
\item[(1)] For any $f,g\in V$, the functions $\max\{f(x),g(x)\}$ and $\min\{f(x),g(x)\}$ belong to $V$.
\item[(2)] For any distinct points $x\neq y$ in $X$, there exists $f\in V$ such that $f(x)\neq f(y)$.
\end{enumerate}
\end{thmm}

Let us go back to $\xcv$. Applying the theorem, we only need to
check that for any distinct points $x, y\in \xcv,$ there exist a
finite extension $E$ of $K_v$, and a continuous function $f$ over
$X_E^{\mathrm{an}}$ such that $f$ takes different values at the
images of $x$ and $y$ in $X_E^{\mathrm{an}}$. This is equivalent
to finding an $E$ such that $x$ and $y$ have different images in
$X_E^{\mathrm{an}}$.

Assume that $x, y\in \xcv$ have the same image in $X_E^{\mathrm{an}}$ for any finite extension $E/K_v$. We are going to show that $x=y$. The problem is local. Assume $\mb(A)$ is an affinoid subdomain of $\xkv$ containing the image of $x$ and $y$. The natural map $\mb(A\widehat\otimes{\CC_v})\rightarrow \mb(A\otimes E)$ is just the restriction of multiplicative semi-norms from $A\widehat\otimes{\CC_v}$ to $A\otimes E$. And thus the semi-norms $x$ and $y$ have the same restriction on $A\otimes E$ for any $E$ by the assumption. But $\bigcup_E A\otimes E=A\otimes{\overline K_v}$ is dense in $A\widehat\otimes{\CC_v}$. It follows that $x$ and $y$ are the same on $A\widehat\otimes{\CC_v}$. That completes the proof.

\subsection{Equidistribution at Infinite Places}
Now we are going to prove Theorem \ref{5.1} and Theorem \ref{5.2} for any archimedean place $v$. We will show Theorem \ref{5.1}, and this is enough by the equivalence relation developed at the end of last subsection. The proof follows the original idea in \cite{SUZ} and \cite{Zh3}, except that we use Theorem \ref{main} to produce small sections instead of the arithmetic Hilbert-Samuel formula.

Assume $v$ is archimedean. Then $\CC_v=\CC$, and $\xcv$ is the usual complex space $X_v(\CC)$. A continuous function $f$ on $\xcv$ is called \textit{smooth} if there is an embedding $\xcv$ in a projective manifold $Y$ such that $f$ can be extended to a smooth function on $Y$. As in \cite{Zh3}, by the Stone-Weierstrass theorem, continuous functions on $\xcv$ can be approximated uniformly by smooth functions.

It suffices to show
$$\lim_{m\rightarrow\infty} \int_{\xcv} f \mu_{x_m}=\frac{1}{\deg_{\lb}(X)}\int_{\xcv} f c_1(\lbb)_v^{n-1}$$
for any smooth real-valued function $f$ on $\xcv$.

For any real function $g$ on $\xcv$ and any metrized line bundle $\mbb=(\mb,\|\cdot\|)$ over $X$, define the \textit{twist} $\mbb(g)=(\mb,\|\cdot\|')$ to be the line bundle $\mb$ over $X$ with the metric $\|s\|'_v=\|s\|_ve^{-g}$ and $\|s\|'_w=\|s\|_w$ for any $w\neq v$. We first prove a lemma.

\begin{lem}\label{5.3}
Assume the above condition, i.e., $\lbb$ is a semipositive metrized line bundle over $X$ and $f$ a smooth real-valued function on $\xcv$. For $\epsilon>0,$ the adelic volume
$$\chi(N\lbb(\epsilon f))
\geq \frac{\chern(\lbb(\epsilon f))^n+O(\epsilon^2)}{n!}N^n+o(N^n),$$
where the error term $O(\epsilon^2)$ is independent of $N$.
\end{lem}
\begin{proof}
See \cite{Zh2} for the definition and basic properties of adelic volumes for adelically metrized line bundles. 

Pick any integral model of $X$. Then $\ob(f)$ is naturally a hermitian line bundle over this model. We can write $\overline\ob(f)=\mbb_1-\mbb_2$ for two ample hermitian line bundles $\mbb_1, \mbb_2$ with positive curvatures over the model. This is the very reason that we assume $f$ is smooth. Still denote  by $\mbb_1, \mbb_2$ the corresponding adelically metrized line bundles over $X$.

We first consider the case that $\lbb$ is induced by a single relatively semipositive model of $(X, \lb)$. We can  assume that the metrics of $\mbb_1, \mbb_2$ and $\lbb$ are induced by line bundles over the same integral model of $X$. This is a standard procedure: any two integral models of $X$ are dominated by the third one, and we pull-back all line bundles to the third one. In this case, we use the notation of each adelic line bundle to denote its corresponding hermitian line bundle. The intersection numbers and adelic volumes are equal to their arithmetic counterpart under this identity. 

Since $\lbb$ is relatively semipositive, there exists a constant $c>0$ such that $\lbb(c)$ is ample. This is a fact we have shown in the part of Notations and Conventions. Then $\lbb(c +\epsilon f)=(\lbb(c)+\epsilon\mbb_1)-\epsilon\mbb_2$ is the difference of two ample hermitian line bundles. Applying Theorem \ref{main}, one gets
\begin{eqnarray*}
\chi_{\sup}(N\lbb(c+\epsilon f))
&\geq&\frac{\chern(\lbb(c)+\epsilon\mbb_1)^n-n\cdot\chern(\lbb(c)+\epsilon\mbb_1)^{n-1}\chern(\epsilon\mbb_2)}{n!}N^n+o(N^n)\\
&=&\frac{\chern(\lbb(c+\epsilon f))^n+O(\epsilon^2)}{n!}N^n+o(N^n)\\
&=&\frac{\chern(\lbb(\epsilon f))^n+cn\deg_{\lb}(X)+O(\epsilon^2)}{n!}N^n+o(N^n).
\end{eqnarray*}
By definition, it is easy to see that
$$\chi_{\sup}(N\lbb(c+\epsilon f))-\chi_{\sup}(N\lbb(\epsilon f))=cN\rank\Gamma(N\lbb)=c\frac{\deg_{\lb}(X)}{(n-1)!}N^n+o(N^n). $$
Thus we obtain the result
$$\chisup(N\lbb(\epsilon f))
\geq \frac{\chern(\lbb(\epsilon f))^n+O(\epsilon^2)}{n!}N^n+o(N^n).$$

Now we come to the general case. It is induced by a limit process. The metric $\lbb$ is a uniform limit of some sequence 
$\{  \|\cdot\|_k \}_k$ of adelic metric on $\lb$, each term of which is induced by a single integral model. Denote $\lbb_k=(\lb, \|\cdot\|_k)$. By uniformity, we can find some constant $c>0$ such that all $\lbb_k(c)$ are ample. We have proved
$$\chi(N\lbb_k(c+\epsilon f)) \geq   \frac{\chern(\lbb_k(c)+\epsilon\mbb_1)^n-n\cdot\chern(\lbb_k(c)+\epsilon\mbb_1)^{n-1}\chern(\epsilon\mbb_2)}{n!}N^n+o(N^n).$$
Set $k\rightarrow \infty$ and we want to get 
$$\chi(N\lbb(c+\epsilon f)) \geq   \frac{\chern(\lbb(c)+\epsilon\mbb_1)^n-n\cdot\chern(\lbb(c)+\epsilon\mbb_1)^{n-1}\chern(\epsilon\mbb_2)}{n!}N^n+o(N^n).$$
Once this is true, the lemma will be proved like the above case.

Dividing both sides of the inequality for $\lbb_k$ by $N^n$, we get
$$\frac{\chi(N\lbb_k(c+\epsilon f))}{N^n} \geq   \frac{\chern(\lbb_k(c)+\epsilon\mbb_1)^n-n\cdot\chern(\lbb_k(c)+\epsilon\mbb_1)^{n-1}\chern(\epsilon\mbb_2)}{n!} +o(1).$$
Now 
$$\lim_{k\rightarrow \infty}  \frac{\chi(N\lbb_k(c+\epsilon f))}{N^n} =\frac{\chi(N\lbb(c+\epsilon f))}{N^n}$$
uniformly. Thus we can find a function $\beta(k)$ with $\lim_{k\rightarrow \infty}\beta(k)=0$ such that 
$$\frac{\chi(N\lbb(c+\epsilon f))}{N^n} \geq \frac{\chi(N\lbb_k(c+\epsilon f))}{N^n} -\beta(k)$$
for $k$ large enough. Hence  
$$\frac{\chi(N\lbb(c+\epsilon f))}{N^n} \geq  -\beta(k)+ \frac{\chern(\lbb_k(c)+\epsilon\mbb_1)^n-n\cdot\chern(\lbb_k(c)+\epsilon\mbb_1)^{n-1}\chern(\epsilon\mbb_2)}{n!} +o(1).$$
Set $k\rightarrow \infty$, and apply Lemma \ref{regularity}. Note that the position of $\epsilon\rightarrow 0$ in the lemma is taken by $k\rightarrow \infty$. We have
$$\frac{\chi(N\lbb(c+\epsilon f))}{N^n} \geq   \frac{\chern(\lbb(c)+\epsilon\mbb_1)^n-n\cdot\chern(\lbb(c)+\epsilon\mbb_1)^{n-1}\chern(\epsilon\mbb_2)}{n!} +o(1).$$
Equivalently, 
$$\chi(N\lbb(c+\epsilon f)) \geq   \frac{\chern(\lbb(c)+\epsilon\mbb_1)^n-n\cdot\chern(\lbb(c)+\epsilon\mbb_1)^{n-1}\chern(\epsilon\mbb_2)}{n!}N^n +o(N^n).$$
Now the argument is the same as the case considered at the beginning of this proof.
\end{proof}

With this lemma, the proof of Theorem \ref{5.1} is the same as the original ones. In fact, fix an archimedean place $w_0$ of $K$. By the adelic Minkowski's theorem (cf. \cite[Appendix C]{BG}), the above lemma implies the existence of a nonzero small section $s\in \Gamma(X, N\lb)$ such that
$$\displaystyle\log\|s\|'_{w_0}\leq
-\frac{\chern(\lbb(\epsilon f))^n+O(\epsilon^2)}{n\deg_{\lb}(X)}N+o(N)=\left(-h_{\lbb(\epsilon f)}(X)+O(\epsilon^2)\right)N+o(N),$$
and $\displaystyle\log\|s\|'_w\leq 0$ for all $w\neq w_0.$ Here $\|\cdot\|'_w$ denotes the metric of $\lbb(\epsilon f)$. Computing the heights of the points in the generic sequence by this section, we get
$$\liminf_{m\rightarrow \infty} h_{\lbb(\epsilon f)}(x_m)\geq h_{\lbb(\epsilon f)}(X)+O(\epsilon^2).$$
By definition,
\begin{eqnarray*}
h_{\lbb(\epsilon f)}(x_m)&=&h_{\lbb}(x_m)+ \epsilon \int_{\xcv} f \mu_{v,x_m},\\
h_{\lbb(\epsilon f)}(X)&=&h_{\lbb}(X)+ \epsilon\frac{1}{\deg_{\lb}(X)} \int_{\xcv} f c_1(\lbb)_v^{n-1}+O(\epsilon^2).
\end{eqnarray*}
Since $$\displaystyle\lim_{m\rightarrow \infty}\height(x_m)=\height(X),$$
we have
$$
\liminf_{m\rightarrow \infty} \int_{\xcv} f \mu_{v,x_m} \geq \frac{1}{\deg_{\lb}(X)}\int_{\xcv} f c_1(\lbb)_v^{n-1}.
$$
Replacing $f$ by $-f$ in the inequality, we get the other direction and thus
$$
\lim_{m\rightarrow \infty}\int_{\xcv} f \mu_{v,x_m}=\frac{1}{\deg_{\lb}(X)}\int_{\xcv} f c_1(\lbb)_v^{n-1}.
$$

\subsection{Equidistribution at Finite Places}
In this subsection, we prove Theorem \ref{5.1} and Theorem \ref{5.2} for any non-archimedean place $v$. We will show Theorem \ref{5.2}, the algebraic case instead of the geometric case. Then Theorem \ref{5.1} is implied by the argument at the end of Section \ref{generic}.

The proof here is parallel to the archimedean case, so the task is to initiate a process which can be run in the same way as in the archimedean case. The key is Gubler's theorem that continuous functions over Berkovich spaces can be approximated by model functions which will be defined later. One can also strengthen Lemma 3.4 and Lemma 3.5 in \cite{Ch2} to prove the result here.

\subsubsection*{Canonical Measures}
The analytic space $\xkv$ is the Berkovich space associated to the
variety $X_{K_v}$ for non-archimedean $v$. The canonical measure
$c_1(\lbb)_v^{n-1}$ is defined by Chambert-Loir \cite{Ch2} using ideas
from the archimedean case. For example, if $\lb_0,\cdots, \lb_d$
are line bundles over $X$ with $v$-adic metrics, and $Z$ is a
closed subvariety of $X$ of dimension $d$, then the local height
formula (for $s_j\in \Gamma(X, \lb_j)$ intersecting properly over
$Z$)
\begin{eqnarray*}
&&(\widehat{\mathrm{div}}(s_0)\cdots\widehat{\mathrm{div}}(s_d)|_Z)_v\\
&=&(\widehat{\mathrm{div}}(s_1)\cdots\widehat{\mathrm{div}}(s_d)|_{\mathrm{div}(s_0|_Z)})_v
-\int_{\xkv}\log \|s_0\|_v c_1(\lb_1)_v\cdots c_1(\lb_d)_v\delta_{Z_{K_v}^{\mathrm{an}}}
\end{eqnarray*}
holds as in the archimedean case. And one also has the global height
$$(\hat{c}_1(\lbb_0)\cdots\hat{c}_1(\lbb_d))|_Z
=\sum_{v}(\widehat{\mathrm{div}}(s_0)\cdots\widehat{\mathrm{div}}(s_d)|_Z)_v,
$$
where the sum is over all places $v$ of $K$.

Denote by $O_{K_v}$ the valuation ring of $K_v$, and by $k_v$ the residue field.
If the $v$-adic metric on $\lbb$ is defined by a single $O_{K_v}$-model $(\mathcal{X}, \lbt)$ with $\mathcal{X}$ normal, then the canonical measure over $\xkv$ has a simple expression
$$c_1(\lbb)_v^{n-1}=\sum_{i=1}^r m_i \deg_{\lbt}(Y_i) \delta_{\eta_i},
$$
where $Y_1, \cdots, Y_r$ are the irreducible components of the special fibre $\xb_{k_v}$, and $m_1,\cdots, m_r$ are their multiplicities, and $\eta_j$ is the unique preimage in $\xkv$ of the generic point of $Y_j$ under the reduction map $\xkv \rightarrow \xb_{k_v}$. Locally, $\eta_j$ is the semi-norm given by the valuation of the local ring of the scheme at the generic point of $Y_j$.

The canonical measures over $\xcv$ have properties similar to the algebraic case.

\subsubsection*{Model Functions}
Let $B$ be a $K_v$-Berkovich space which is Hausdorff, compact and strictly $K_v$-analytic. There is a notion of \textit{formal $O_{K_v}$-model} for $B$, which is an \textit{admissible formal $O_{K_v}$-scheme} with \textit{generic fibre} $B$. For the basics of formal models we refer to \cite{Ra} and \cite{BL}. Let $M$ be a line bundle over $B$. Among the $K_v$-metrics over $M$, there are some called \textit{formal metrics} by Gubler \cite{Gu}. They are induced by formal models of $(B, M)$.

\begin{definition}\label{5.4}
A continuous function over $B$ is called a model function if it is equal to $-\log\|1\|^{1/l}$ for some nonzero integer $l$ and some formal metric $\|\cdot\|$ over the trivial bundle of $B$.
\end{definition}

It is easy to see that all model functions form a vector space. The following theorem is due to Gubler \cite[Theorem 7.12]{Gu}.

\begin{thmm}[Gubler]
The vector space of model functions on $B$ is uniformly dense in the ring of real-valued continuous functions on $B$.
\end{thmm}

Now let's come back to our situation: $X$ is a projective space over $K$ and $\xkv$ is the corresponding Berkovich space at $v$. To compute heights, we work on \textit{global projective $O_K$-model} of $(X, \ob_X)$ in the usual sense, i.e. a pair $(\xb, \mb)$ consisting of an integral scheme $\xb$ projective and flat over $O_K$ with generic fibre $X$, and a line bundle $\mb$ over $\xb$ which extends $\ob_X$.

A global projective $O_K$-model gives a formal $O_{K_v}$-model by completion with respect to the ideal sheaf $(\varpi)$ where $\varpi$ is a uniformizer of $O_{K_v}$. Thus it induces a formal metric over $\ob_X$, which is compatible with the adelic metric defined by Zhang. Now we are going to show that all formal metrics arise in this way.

\begin{lem}\label{5.5}
All formal metrics over the trivial bundle of $\xkv$ are induced by global projective $O_K$-models. Thus all model functions are induced by global projective $O_K$-models.
\end{lem}
\begin{proof}
Let $(\xb, \mb)$ be any formal $O_{K_v}$-model. Fix a global projective $O_K$-model $\xb_0$ of $X$ and denote by
$\widehat\xb_0$ its completion at $v$. Then $\widehat\xb_0$ gives another formal $O_{K_v}$-model of $\xkv$. By Raynaud's result (cf. \cite[Theorem 4.1]{BL}) on the category of formal models, there exist two admissible formal blowing-ups
$\xb' \rightarrow \widehat\xb_0$ and $\phi: \xb' \rightarrow \xb$, both of which induce isomorphisms over $\xkv$. Then $(\xb', \mb')$ induces the same formal metric as $(\xb, \mb)$, where $\mb'=\phi^*\mb$.

Denote by $\ib$ the coherent ideal sheaf for the blowing-up $\xb' \rightarrow \widehat\xb_0$. By the formal GAGA (cf. EGA III.1, Section 5), $\ib$ comes from a coherent ideal sheaf of the projective $O_{K_v}$-variety $(\xb_0)_{O_{K_v}}$. We still denote it by $\ib$. We can also consider $\ib$ as a coherent ideal sheaf of $\xb_0$, since $\ib$ contains some power of the maximal ideal $(\varpi)$ of $O_{K_v}$. Let $\xb''$ be the blowing-up of $\xb_0$ with respect to $\ib$. Then $\xb''$ is a global projective $O_K$-model, and the completion at $v$ of $\xb''$ gives $\xb'$.

Now it remains to find a model of $\mb'$ over $\xb''$. One can descend $\mb'$ to a line bundle over $\xb''_{O_{K_v}}$ by formal GAGA, and we still denote it by $\mb'$. Let $D$ be a divisor on $\xb''_{O_{K_v}}$ defined by any rational section of $\mb'$. Since $\mb'$ is trivial over the generic fibre of $\xb''_{O_{K_v}}$, there exists a positive integer $r$ such that $r \ \mathrm{div}(\varpi)+D$ is effective, where $\mathrm{div}(\varpi)$ is the whole special fibre. Replacing $D$ by $r\ \mathrm{div}(\varpi)+D$, we assume that $D$ is effective.

Let $\jb$ be the ideal sheaf of $D$ in $\xb''_{O_{K_v}}$. Then $\jb$ is invertible over $\xb''_{O_{K_v}}$. We can also considere $\jb$ as a coherent ideal sheaf of $\xb''$. If $\jb$ is invertible over $\xb''$, then $(\xb'', \jb^{\otimes (-1)})$ is a desired global projective model which gives the same metric as $(\xb', \mb')$ does. Otherwise, consider the blowing-up $\pi: \xb'''\rightarrow\xb''$ with respect to $\jb$. Then $\pi^{-1}\jb$ is invertible over $\xb'''$, and $(\xb''', (\pi^{-1}\jb)^{\otimes (-1)})$ gives what we want. (In fact, one can show that $\jb$ is invertible over $\xb''$ by this blowing-up.)
\end{proof}

\begin{remark}
It is possible to work directly on global projective $O_K$-models and show that the model functions defined by them are uniformly dense, which will be enough for our application. Of course, it still follows Gubler's idea in proving the density theorem. Use the Stone-Weierstrass theorem. Pick any initial projective model, blow-up it suitably to get separation of points, and use certain combinatorics and blowing-ups to prove the model functions are stable under taking maximum and minimum.
\end{remark}

\subsubsection*{A Description of the Berkovich Space}
Using model functions, Zhang \cite{Zh5} constructed the Berkovich space $\xkv$ in an elementary way. For any projective variety $X$ over $K$, let $V$ be the vector space of all model functions coming from (varying) projective $O_{K_v}$-models of $X$. Each element of $V$ is considered as a map from $|X_{K_v}|$ to $\RR$, where $|X_{K_v}|$ is purely the underlying space of the scheme. Now take $R(X_{K_v})$ to be the completion under the supremum norm of the ring generated by $V$. Then we have
$$\xkv=\mathrm{Hom}(R(X_{K_v}), \RR),$$
where Hom is taking all continuous homomorphisms.

In fact, by the density of model functions, $R(X_{K_v})$ is exactly the ring of continuous functions over the compact Hausdorff space $\xkv$. Therefore its spectrum recovers $\xkv$. The same construction is valid for $\xcv$.

\subsubsection*{Proof of Equidistribution}
Now we are ready to prove Theorem \ref{5.2} when $v$ is non-archimedean. By the density theorem proved above, it suffices to show
$$\lim_{m\rightarrow\infty} \int_{\xkv} f\mu_{v, x_m}=\frac{1}{\deg_{\lb}(X)}\int_{\xkv} f c_1(\lbb)_v^{n-1}$$
for any model function $f=-\log \|1\|_v$ induced by a projective $O_K$-model $(\xb, \mb)$.

Denote by $\ob(f)$ the trivial line bundle $\ob_X$ with the adelic metric given by the model $(\xb, \mb)$, i.e., the metric such that $\|1\|_v=e^{-f}$ and $\|1\|_w=1$ for any $w\neq v.$ Define the \textit{twist} $\lbb(\epsilon f)=\lbb+\epsilon \ob(f)$ for any positive rational number $\epsilon$. Note that we even have exactly the same notation as in the archimedean case.

Over $\xb$, the line bundle $\mb$ is a difference of two ample hermitian line bundles. It follows that $\ob(f)$ is a difference of two ample metrized line bundles. This tells why we spend so much energy proving the density of model functions induced by global models.

Now everything including Lemma \ref{5.3} follows exactly in the same way. In particular, we have
$$\liminf_{m\rightarrow \infty} h_{\lbb(\epsilon f)}(x_m)\geq h_{\lbb(\epsilon f)}(X)+O(\epsilon^2).$$
By the definition of our metrics and intersections,
\begin{eqnarray*}
h_{\lbb(\epsilon f)}(x_m)&=&h_{\lbb}(x_m)+ \epsilon \int_{\xkv} f \mu_{v,x_m},\\
h_{\lbb(\epsilon f)}(X)&=&h_{\lbb}(X)+ \epsilon \frac{1}{\deg_{\lb}(X)} \int_{\xkv} f c_1(\lbb)_v^{n-1}+O(\epsilon^2).
\end{eqnarray*}
The variational principle follows exactly in the same way.

\subsection{Equidistribution of Small Subvarieties}\label{subvarieties}
Szpiro-Ullmo-Zhang's equidistribution theorem was generalized to equidistribution of small subvarieties by Baker-Ih \cite{BI} and Autissier \cite{Au2}. Now we will generalize our theory to small subvarieties in the same manner. We omit  the proof, since it follows the treatment of \cite{Au2} by using the variational principle. We will only formulate the result in the geometric case, though it is immediate for both cases.

Suppose we are in the situation of Theorem \ref{5.1}. More precisely, $X$ is a projective variety of dimension $n-1$ over a number field $K$, and $\lbb$ is a metrized line bundle over $X$ such that $\lb$ is ample and the metric is semipositive.

By a \textit{subvariety} of $X$, we mean a reduced and irreducible closed subscheme defined over $\overline K$. For any subvariety $Y$ of $X$, define its \textit{height} to be
$$\height(Y)=\frac{\chern(\lbb)^{\dim Y+1}|_{\overline Y}}{(\dim Y+1)\deg_{\lb}(\overline Y)},$$
where $\overline Y$ is the closure of $Y$ in the scheme $X$.

Then $\overline Y_{\overline K}$ splits into a finite set of subvarieties in $X_{\overline K}$. We denote this set by $O(Y)$, and call it the \textit{Galois orbit} of $Y$. For any $Z\in O(Y)$, the associated analytic space $Z_{\CC_v}^{\mathrm{an}}$ is a closed subspace of $\xcv$. Thus we can also view $O(Y)$ as a finite set of closed analytic subspace of $\xcv$ for any place $v$.

Now define \textit{the probability measure associated to $Y$} by
$$\displaystyle\mu_{v,Y}=\frac{1}{\deg_{\lb}(\overline Y)}\sum_{Z\in O(Y)}\delta_{Z_{\CC_v}^{\mathrm{an}}}c_1(\lbb|_Z)_v^{\dim Y},$$
where $c_1(\lbb|_Z)_v^{\dim Y}$ is the $v$-adic canonical measure over $Z_{\CC_v}^{\mathrm{an}}$,  and $\delta_{Z_{\CC_v}^{\mathrm{an}}}c_1(\lbb|_Z)_v^{\dim Y}$ sends a continuous function $f: \xcv\rightarrow \CC$ to
$\displaystyle\int_{Z_{\CC_v}^{\mathrm{an}}} f c_1(\lbb|_Z)_v^{\dim Y}.$

We need an additional assumption: $h_{\lbb}(Y)\geq h_{\lbb}(X)$ for any subvariety $Y$ of $X$. We will see later that for dynamical systems $h_{\lbb}(X)=0$ and $h_{\lbb}(Y)\geq 0$ is always true. If $\lbb$ is an ample metrized line bundle, the assumption is equivalent to $h_{\lbb}(x)\geq h_{\lbb}(X)$ for any point $x$ of $X$ by the successive minima of Zhang \cite[Theorem 1.10]{Zh2}.

With the same notions of small sequences, generic sequences and equidistribution as in Section \ref{generic}, we have:

\begin{thm}[Equidistribution of Small Subvarieties]\label{5.6}
Suppose $X$ is a projective variety of dimension $n-1$ over a number field $K$, and $\lbb$ is a metrized line bundle over $X$ such that $\lb$ is ample and the metric is semipositive. Assume $h_{\lbb}(Y)\geq h_{\lbb}(X)$ for any subvariety $Y$ of $X$.  Let $\{Y_m\}$ be an infinite sequence of subvarieties of $X$ which is generic and small. Then for any place $v$ of $K$, the Galois orbits of the sequence $\{Y_m\}$ are equidistributed in the analytic space $\xcv$ with respect to the canonical measure $d\mu_v=c_1(\lbb)_v^{n-1}/\deg_{\lb}(X)$.
\end{thm}

\subsection{Equidistribution over Algebraic Dynamics}
\label{dynamics}
The equidistribution theorems treated in previous subsections have direct consequences in algebraic dynamics. For a complete introduction to the basics and equidistribution of algebraic dynamics, we refer to \cite{Zh5}. And we will only state the equidistribution of small points in the geometric case.

Let $K$ be a number field. Let $X$ be a projective variety over $K$, and $\phi: X \rightarrow X$ be a morphism polarized by an ample line bundle $\lb$ over $X$, meaning that $\phi^*\lb \cong \lb^{\otimes q}$ for some integer $q>1$. Then $(X, \phi, \lb)$ is called an \textit{algebraic dynamical system}.

Fix an isomorphism $\alpha: \phi^*\lb \cong \lb^{\otimes q}$. By \cite{Zh2}, there exists a unique semipositive metric over $\lb$ which makes $\alpha$ an isometry. Actually, it can be obtained by Tate's limit like the canonical height and the canonical measure in Section 1. This metric is called \textit{the canonical metric}. Denote by $\lbb$ the line bundle $\lb$ endowed with this metric. For any place $v$ of $K$, one has \textit{the canonical measure} $c_1(\lbb)^{n-1}_v$  and \textit{the canonical probability measure} $d\mu_{v,\phi}:=c_1(\lbb)_v^{n-1}/\deg_{\lb}(X)$ over $\xcv$.

Using the canonical metric, we define the \textit{canonical height} of a subvariety $Y$ by
$$\cheight(Y)=\height(Y)=\frac{\chern(\lbb)^{\dim Y+1}|_{\overline Y}}{(\dim Y+1)\deg_{\lb}(\overline Y)},$$
as in the previous subsection. It is the same as the one defined by Tate's limit.

Now we use the same notions of small sequences, generic sequences and equidistribution as in Section \ref{generic}. A closed subvariety $Y$ of $X$ is called \textit{preperiodic} if the orbit $\{Y, \phi(Y), \phi^2(Y), \cdots\}$ is finite. Note that $\cheight(X)=\height(X)=0$ since $X$ is preperiodic, and thus a small sequence really has heights going to zero. Now the following theorem is just a dynamical version of Theorem \ref{5.1}.

\begin{thm}[Dynamical Equidistribution of Small Points]\label{5.7}
Let $(X, \phi, \lb)$ be an algebraic dynamical system over a number field $K$, and $\{x_m\}$ be an infinite sequence of
algebraic points of $X$ which is generic and small. Then for any place $v$ of $K$, the Galois orbits of the sequence $\{x_m\}$ are equidistributed in the analytic space $\xcv$ with respect to the canonical probability measure $d\mu_{v, \phi}=c_1(\lbb)_v^{n-1}/\deg_{\lb}(X)$.
\end{thm}

\begin{remark}
Following the formulation in Section \ref{subvarieties}, we have equidistribution of small subvarieties over a dynamical system. 
\end{remark}

As in \cite{SUZ}, this result gives the equivalence between the dynamical Bogomolov conjecture and the strict equidistribution of small points.

\begin{conj}[Dynamical Bogomolov Conjecture]
Let $Y$ be an irreducible closed subvariety of $X$ which is not preperiodic. Then there exists a positive number $\epsilon>0,$ such that the set $\{x\in Y(\overline K): \cheight(x)<\epsilon\}$ is not Zariski dense in $Y$.
\end{conj}

\begin{remark}
The known cases of this conjecture are: the case of multiplicative
groups by Zhang \cite{Zh1}, the case of abelian varieties proved by
Ullmo \cite{Ul} and Zhang \cite{Zh3}, and the almost split semi-abelian case
proved by Chambert-Loir \cite{Ch1}. The general case without any group
structure is widely open.
\end{remark}

A sequence $\{x_m\}_{m\geq 1}$ of algebraic points in $X$ is call \textit{strict} if no infinite subsequence of $\{x_m\}$ is contained in a proper preperiodic subvariety of $X$. The strict equidistribution is the following:

\begin{conj}[Dynamical Strict Equidistribution of Small Points]
Let $\{x_m\}$ be an infinite sequence of algebraic points of $X$ which is strict and small. Then for any place $v$ of $K$, the Galois orbits of the sequence $\{x_m\}$ are equidistributed in the analytic space $\xcv$ with respect to the canonical probability measure $d\mu_{v, \phi}=c_1(\lbb)_v^{n-1}/\deg_{\lb}(X)$.
\end{conj}

\begin{cor}\label{5.8}
Dynamical Bogomolov Conjecture $\Longleftrightarrow$ Dynamical Strict Equidistribution of Small Points.
\end{cor}

By a result of Bedford-Taylor \cite{BT} and Demailly \cite{De}, the support of the canonical measure is Zariski dense in $Y_{\CC_v}^{\mathrm{an}}$ for archimedean $v$. See also \cite[Theorem 3.1.6]{Zh5} for example. With this result, the corollary is easily implied by Theorem \ref{5.7}.

Address: \textit{Department of Mathematics, Columbia University,
New York, NY 10027.}

Email: \textit{yxy@math.columbia.edu.}
\end{document}